\documentclass[%
a4paper,%
defaults,%
arxiv,%
]%
{myclass}

\usepackage{nicefrac}
\usepackage{booktabs}
\usepackage{tikzexternal}

\usepackage[%
  bb,%
  ]{mymacros}

\usepackage[nowarnings]{hopf}
\usepackage{hopf.abbrev}

\tikzsetexternalprefix{free-}
\tikzsetfigurename{figure}

\newtheorem{examples}[theorem]{Examples}

\newcommand{\eHom}[3]{\underline{#1}(#2, #3)}

\let\swprod\boxdot

\makeatletter
\renewcommand{\textvCStobj}{%
  \tCo{}\tSo{}\hopfLower{M}onoid\hopfSingle{}{s}\hopf@matharg{\hopf@arg}
 in \(\vcat\)}

\renewcommand{\textgobj}{%
  \hopfLower{G}roup\hopfSingle{}{s}\hopf@matharg{\hopf@arg}}

\renewcommand{\textaobj}{%
  \hopfLower{A}belian \hopfLower{G}roup\hopfSingle{}{s}\hopf@matharg{\hopf@arg}}

\renewcommand{\textrobj}{%
\hopfLower{R}ing\hopfSingle{}{s}\hopf@matharg{\hopf@arg}}

\renewcommand{\textbobj}{%
\hopfLower{C}ommutative \(\robj\)\enhyp{}\hopfLower{A}lgebra\hopfSingle{}{s}\hopf@matharg{\hopf@arg}}

\renewcommand{\textsobj}{%
\hopfLower{S}et\hopfSingle{}{s}\hopf@matharg{\hopf@arg}}

\makeatother

\mybibliographystyle

\mysubjclass[55S25,16W30]{18D35}

\mykeywords{unstable cohomology operations, Hopf rings}

\mytitle{%
  Tall\enhyp{}Wraith Monoids}
\mydate{\today}

\myauthor{Andrew Stacey}
\myaddress{Institutt for Matematiske Fag\\NTNU\\7491 Trondheim\\Norway}
\myemail{stacey@math.ntnu.no}
\myurl{http://www.math.ntnu.no/\(\sim\)stacey}

\myauthor{Sarah Whitehouse}
\myaddress{School of Mathematics and Statistics\\University of Sheffield\\Sheffield\\UK}
\myemail{S.Whitehouse@sheffield.ac.uk}
\myurl{http://www.sarah-whitehouse.shef.ac.uk}

\mythanks{The authors acknowledge the support of the EPSRC, grant no.: GR/S76823
/01.}

\mykeywords{Unstable cohomology operations -- Hopf ring.}

\myabstract{%
  Tall\enhyp{}Wraith monoids were introduced in \cite{MR2559638} to describe the algebraic structure on the set of unstable operations of a suitable generalised cohomology theory.
In this paper we begin the study of Tall\enhyp{}Wraith monoids in an algebraic and categorical setting.
We show that for \(\vcat\) a variety of algebras, applying the \free{v} to a monoid in \(\scat\) produces a Tall\enhyp{}Wraith monoid.
We also study the example of the Tall\enhyp{}Wraith monoid defined by the self set\hyp{}maps of a finite ring, an example closely related to the original motivation for Tall\enhyp{}Wraith monoids.
}

\begin{document}

\mymaketitle

\section{Introduction}
\label{sec:intro}

In \cite{MR2559638}, we introduced the notion of a \emph{Tall\enhyp{}Wraith monoid} as a way of describing the algebraic structure on the set of unstable operations of a suitable generalised cohomology theory.
In brief, for a variety of algebras, \(\vcat\), a \emph{Tall\enhyp{}Wraith \(\vcat\)\enhyp{}monoid} is ``that which can act on \valgs''.
The connection to cohomology theories is almost immediate since the cohomology of a space is an algebra over the coefficient ring.
We say ``almost immediate'' since one has to take into account two additional pieces of structure: the grading and the filtration.
Doing this comprises the main work of \cite{MR2559638}.
Our purpose in laying this foundation is to give direct descriptions, with  a ``generator\enhyp{}relation'' flavour, of the structure on the sets of unstable operations for suitable cohomology theories.
What led us to consider this was the work of \cite{assw4} where we related stable and unstable operations of suitable mod \(p\) cohomology theories.
We sought a simple setting for this result but although this exists for stable operations we found we had to invent it ourselves for unstable operations. 

A simple example of a Tall\enhyp{}Wraith monoid is that of a ring.
The category of abelian groups forms a variety of algebras and ``that which can act on abelian groups'' is simply a ring; unital, but not necessarily commutative.
However, this example hides some of the complexity of the theory.
Because the category of abelian groups is particularly nice, the associated category of Tall\enhyp{}Wraith monoids is again a variety of algebras.
This is not true in general and this complicates matters somewhat.

A ``generators and relations'' description of, say, a ring can be understood in categorical terms as writing that ring as the coequaliser of two free rings.
That this works relies heavily on the properties of the ``free ring'' functor from sets to rings.
Whilst this works well for any variety of algebras, the categories of Tall\enhyp{}Wraith monoids are not, in general, varieties of algebras so this cannot be assumed to work for those.

A more general view of ``generators and relations'' descriptions is that they are simply a way of describing the unknown object in terms of certain standard, or reference, objects.
This view, being more general, has a better chance of succeeding for Tall\enhyp{}Wraith monoids.
For it to work, we need a good source of reference objects.
The aim of this paper is to provide such a source and to begin the study of Tall\enhyp{}Wraith monoids in an algebraic and categorical setting.

\medskip

We examine three sources of Tall\enhyp{}Wraith monoids: using free functors from \(\scat\), free functors where the base is a commutative variety of algebras, and something that we call \emph{toy cohomology theories}: namely, cohomology theories without the interesting parts  (filtration and grading).
The first two of these are what we mean by a ``good source of reference objects''.
These produce Tall\enhyp{}Wraith monoids with a particularly simple structure which can be easily analysed using already developed tools.
The last, toy cohomology theories, approaches the study from the other end.
Our intention in including it is to give an example of relating a Tall\enhyp{}Wraith monoid found ``in the wild'' to one from our ``tame'' collection.

\medskip

Let us now describe our constructions in a little more detail.
Our main tool for building these reference objects is, in fact, the free functor for the corresponding variety of algebras.
To describe how this works, we need to be more precise in our definitions.

\begin{defn}
Let \(\vcat\) be a variety of algebras.
A \emph{Tall\enhyp{}Wraith \(\vcat\)\enhyp{}monoid} is a monoid in \(\vCvcat\), \vCvcat.
\end{defn}

Of course, to show that this is well\hyp{}defined we need to give \vCvcat a monoidal structure.
The proof of this was omitted from \cite{MR2559638} and we include it here.

Our first result shows that the \free{v} lifts from \(\vcat\) to \(\vCvcat\) and is monoidal.

\begin{thm}
 \label{th:twmon}
Let \vcat be a variety of algebras.
The \free{v}, \(\free{v} \colon \scat \to \vcat\), has a lift to a strong monoidal functor \(\free[hat]{v} \colon \scat \to \vCvcat\).
\end{thm}

Thus we have our first source of Tall\enhyp{}Wraith \(\vcat\)\enhyp{}monoids.

\begin{cor}
\label{cor:freemon}
If \(M\) is a monoid then \(\free[hat]{v}(M)\) is a Tall\enhyp{}Wraith \(\vcat\)\enhyp{}monoid.
\end{cor}

The proof of this theorem is straightforward and uses no difficult techniques.
Indeed, we have seen parts of this theorem implicitly used in several instances but have not seen an actual proof.
Of particular note is the paper \cite{pf2} which contains this result  applied to the monoid with one element.

\medskip

There is an analogue of this result for \emph{graded} algebras; these are sometimes called \emph{many\hyp{}sorted algebras} in the literature.
A graded algebra has components corresponding to some (fixed) grading set, say \(\iset\), and the inputs of the operations are constrained to come from certain components.
We have a graded version of Theorem~\ref{th:twmon}.

\begin{thm}
 \label{th:twmongrade}
 Let \(\Gvcat\) be a variety of graded algebras.
 \GvCGvcatu has a monoidal structure,
  \((\GvCGvcat, \twprod, I)\).

 Let
  \(\free[mZo]{Gv} \colon \ZZscat \to \ZGvcat\)
 be the obvious extension of the \free{Gv}.
 Let \(\diag \colon \scat \to \ZZscat\) be the diagonal functor.
 Then the functor
  \(\free[mZo]{Gv}\diag \colon \scat \to \ZGvcat\),
 has a lift to a strong monoidal functor
  \(\free[hat]{Gv} \colon \scat \to \GvCGvcat\).
\end{thm}

\medskip

The second source of Tall\enhyp{}Wraith monoids that we consider comes from considering two varieties of algebras, say \(\vcat\) and \(\ucat\), with the property that every \ualg is again a \valg.
That is, that there is a forgetful functor \(\ucat \to \vcat\).
There is, then, a left adjoint \(\free{uv} \colon \vcat \to \ucat\).
Moreover, we assume that \(\vcat\) is what is sometimes called a \emph{commutative} variety of algebras, such as modules over a commutative ring.
This is enough to ensure that the category \(\vcat\) is closed symmetric monoidal.
This monoidal structure interacts well with the Tall\enhyp{}Wraith monoidal structure on \(\uCucat\) and so the free functor \(\free{uv}\) provides another source of Tall\enhyp{}Wraith monoids.

\medskip

Finally, we consider toy cohomology theories.
By removing the grading and the filtration from the idea of a cohomology theory, we are left with a representable contravariant functor from \(\scat\) to \(\rcat\).
The representing object, say \(\robj\), is thus a ring.
If \(\robj\) is finite, the set of ``operations'' of this theory, \(\Hom{\scat}{\abs{\robj}}{\abs{\robj}}\), is a Tall\enhyp{}Wraith monoid in \bcat.
We can relate this to the free \bobj on a single generator and so describe \(\Hom{\scat}{\abs{\robj}}{\abs{\robj}}\) using \(\robj\lb \selt \rb\).

\medskip

As we said above, we have seen implicit uses of some of the work of this paper in the literature.
Listing them all would be impossible, but two stand out as being worth mentioning because the main themes of those papers are closely related to this one.
Those papers are \cite{dtgw} and \cite{jbbw}.

The first of these is the earliest trace that we have been able to find of these ideas and thus the source of the name \emph{Tall\enhyp{}Wraith monoid}.
(Of course, the authors of \cite{dtgw} did not give them that name!)
In that article, the authors introduced the notion of a \emph{biring} and studied its structure and the structure of the corresponding category.
In the modern parlance of general or universal algebra, a biring is a \emph{\co{ring} object in the category of rings}; here, ``ring'' means commutative, unital ring.
The authors of \cite{dtgw} also consider monoids in the category of birings which they call \emph{biring triples}.
They show that the free ring on a group (a.k.a.\ the group ring) is a biring triple.
This is a specific example of Theorem~\ref{th:twmon}.

All of this, including the example of group rings, was generalised in \cite{jbbw} from the category of rings to that of commutative algebras over a ring; the authors of \cite{jbbw} introduced the name \emph{plethory} for the generalisation of a biring triple.

\medskip

This paper is organised as follows.
Section~\ref{sec:genalg} provides the background on general algebra that we shall need in the rest of the paper.
In Section~\ref{sec:twmon} we establish our starting point by showing that \(\vCvcat\) (\vCvcat) is a monoidal category and giving the definition of a Tall\enhyp{}Wraith \(\vcat\)\enhyp{}monoid.
In Section~\ref{sec:free} we prove Theorems~\ref{th:twmon} and~\ref{th:twmongrade}.
Section~\ref{sec:abelian} concerns the source of examples of \vCvobjs using free functors starting from a commutative variety of algebras.
In Section~\ref{sec:toys} we describe how to obtain a Tall\enhyp{}Wraith monoid from a toy cohomology theory.
In Section~\ref{sec:examples}, we carry out these constructions in  common situations and conclude by comparing the Tall\enhyp{}Wraith monoids from certain toy cohomology theories with the ones from the free functors.
Lastly, in Appendix~\ref{ap:notation} we have a table of the various categories with a brief explanation of the notations that we employ.

\section{General Algebra}
\label{sec:genalg}

In this section we record the main results of general algebra that we need.
The results quoted in this section are all standard results from that field.
For ungraded algebra objects in \scat, these results can be found in any good introduction to the subject, for example \cite{gb}.
The more general cases can be found in the wider literature, for example in \cite{wkjmmpis}.
We record these results here, mostly without proofs, to establish notation and as a quick reference for the rest of the paper.

There are many ways to think of the ideas in general algebra, from Lawvere theories to operations and identities.
Our overall goal is to understand the structure of the operations on a  cohomology theory with a view to being able to study particular operations.
Thus we use language from the ``operations and identities'' view.

We note that we are considering finitary algebraic theories \emhyp{} that is, Lawvere theories \emhyp{} and that as we wish to focus on the corresponding ``sets with structure'', we shall most often refer to the associated ``variety of algebras'' by which we mean the category of models of the theory in \scat.

We start by summarising the results that we need in the arena of ungraded algebras, also known as \emph{single\hyp{}sorted} or \emph{homogeneous} algebras.
We then explain how this generalises to graded algebras, also known as \emph{many\hyp{}sorted} or \emph{heterogeneous} algebras.

We use the exponential notation to denote the iterated product of an object in a category, assuming that these products exist, so \(Y^n\) is the \(n\)\hyp{}fold product of \(Y\) and, for a \sobj[\sobj] \(Y^{\sobj}\) is the product of copies of \(Y\) indexed by \(\sobj\); in particular \(Y^0\) and \(Y^\emptyset\) both denote the terminal object in the category (assuming that it exists).

\subsection{Ungraded Algebra}

\begin{defn}
A \emph{Lawvere theory}, \vcat, is a category with finite cartesian products in which every object is isomorphic to a finite cartesian power of a certain distinguished object.

A \emph{model} of a Lawvere theory, \vcat, in a category, say \dcat, with finite products is a product\enhyp{}preserving functor \(\vcat \to \dcat\).

Clearly, one then obtains the category of models of \(\vcat\) in \(\dcat\) which we denote by \(\dvcat\).
\end{defn}

Thus the key data in a Lawvere theory are the morphisms from a general cartesian power of the distinguished object to that object; these are called the \emph{operations} of the theory.
In a presentation of the Lawvere theory, one starts with a set of operations (called \emph{primary} operations) which, by composition, generate all of the operations.
Wherever two distinct sequences of compositions of primary operations produce the same operation, one has an \emph{identity}.

A Lawvere theory is determined by its category of models in \scat, for which we use the following terminology.

\begin{defn}
A \emph{variety of algebras} is the category of models in \scat of a Lawvere theory.
\end{defn}

We shall henceforth identify the Lawvere theory with its variety of algebras and simply state ``let \vcat be a variety of algebras''.
It is also possible when considering a model of a Lawvere theory to focus attention on the image of the distinguished object rather than on the image of the whole functor.
This leads one to think of \(\dvcat\), the category of models of the Lawvere theory, as \dvcat and to talk of a \dvobj rather than of a model of \vcat.
To completely remove any mention of the underlying Lawvere theory from the discussion, we use the following characterisation of \dvobjs.

\begin{proposition}
 \label{prop:conliftid}
 To give an \dobj[\abs{\doobj}] the structure of a \valgobj is equivalent to giving a lift of the contravariant hom\hyp{}functor
  \(\Hom{\dcat}{-}{\abs{\doobj}} \colon \dcat \to \scat\)
 to a functor \(\dcat \to \vcat\).
 \noproof
\end{proposition}

For a \dvobj[\dvobj] and \dobj[\dobj] we shall write
\[
  \Hom{\dcat}{\dobj}{\dvobj}
\]
for the \valg with underlying \sobj
\[
  \Hom{\dcat}{\dobj}{\abs{\dvobj}}
\]
where \(\abs{\dvobj}\) is the underlying \dobj of \(\dvobj\).

The main work in this paper involves establishing adjunctions between various functors.
The starting point is knowing that the obvious forgetful functor, \(\dvcat \to \dcat\), has a left adjoint.
The following result lists simple conditions that guarantee its existence.

\begin{proposition}
 \label{prop:vfree}
 Let \(\dcat\) be an (\(\m{E}\), \(\m{M}\)) category for some classes \(\m{E}\) of epimorphisms and \(\m{M}\) of monomorphisms.
 Suppose that:
\begin{enumerate}
\item \(\dcat\) has finite products,
\item \(\dcat\) is closed under filtered \co{}limits,
\item finite products commute with filtered \co{}limits,
\item \(\m{E}\) is closed under finite products, and
\item \(\dcat\) is \(\m{E}\)\hyp{}\co{well}\hyp{}powered.
\end{enumerate}

 Then the forgetful functor \(\dvcat \to \dcat\) has a left adjoint,
  \(\free{v} \colon \dcat \to \dvcat\),
 which is called the \free{v}.
 \noproof
\end{proposition}

Dual to \valgobjs are \Cvalgobjs.

\begin{defn}
 Let \(\vcat\) be a variety of algebras, \(\dcat\) a category with finite \co{}products.
 A \emph{\dCvobj{}} is a \Odvobj.
 A \emph{morphism of \dCvobjs} is a morphism in \(\dcat\) which intertwines the \dCvobj structures.
 We denote \dCvcat by \(\dCvcat\).
\end{defn}

The morphisms are chosen such that there is an isomorphism of categories
\[
 \dCvcat \cong \mOo{(\Odvcat)}
\]
and there is an obvious covariant forgetful functor \(\dCvcat \to \dcat\).
The analogue of Proposition~\ref{prop:conliftid} is the following.


\begin{proposition}
 \label{prop:covlift}
 To give an \dobj[\abs{\dCvobj}] the structure of a \Cvalgobj is equivalent to giving a lift of the covariant hom\hyp{}functor
  \(\Hom{\dcat}{\abs{\dCvobj}}{-} \colon \dcat \to \scat\)
 to a functor \(\dcat \to \vcat\).
 \noproof
\end{proposition}

Analogously to earlier,
 \(\Hom{\dcat}{\dCvobj}{\dobj}\)
will denote the \valg with underlying \sobj
 \(\Hom{\dcat}{\abs{\dCvobj}}{\dobj}\).

\subsection{Graded Algebras}

We turn now to graded algebras.
A graded algebra has components indexed by some (fixed) set and its operations go from components to components rather than being globally defined.
The theory of graded algebras is very similar to that of ungraded algebras in that all of the results for ungraded algebras have their counterparts in graded algebras.
The thing to keep track of is where the indexing set plays a part.
In particular, a graded variety of algebras does not have an underlying \sobj but rather an underlying family of \sobjs.

Thus to discuss graded algebras in some category we first need to establish notation for graded objects.
Let us fix some (non\hyp{}empty) grading set \(\iset\).
We shall regard this both as a set and as a (small) discrete category.
We write \(\Zdcat\) for \Zdcat.
As \(\iset\) is a discrete category, there is no distinction between covariant and contravariant functors from \(\iset\).

An \Zdobj[\Zdobj] represents both a covariant and a contravariant functor \(\dcat \to \Zscat\) via
\begin{align*}
 \cov{\Zdobj}(\dobj') &%
 = \big(\isetelt \mapsto \Hom{\dcat}{\Zdobj(\isetelt)}{\dobj'}\big), \\
 \con{\Zdobj}(\dobj') &%
 = \big(\isetelt \mapsto \Hom{\dcat}{\dobj'}{\Zdobj(\isetelt)}\big).
\end{align*}

As the distinction between graded and ungraded varieties of algebras is something that we wish to keep clear, we shall use write graded varieties of algebras with a superscript star, as in \(\Gvcat\).
It is to be emphasised that \(\Gvcat\) and \(\vcat\) bear no relation to each other, and we shall never use both \(\Gvcat\) and \(\vcat\) in the same context.

In the graded case, Proposition~\ref{prop:conliftid} becomes the following result.

\begin{proposition}
 \label{prop:conliftidgrade}
 To give an \Zdobj[\abs{\dGoobj}] the structure of a \Gvalgobj is equivalent to giving a lift of the contravariant hom\hyp{}functor
  \(\Hom{\dcat}{-}{\abs{\dGoobj}} \colon \dcat \to \Zscat\)
 to a functor \(\dcat \to \Gvcat\).
 \noproof
\end{proposition}

As before, if \(\dGvobj\) is a \dGvobj and \(\dobj\) is an \dobj we  write
\[
 \Hom{\dcat}{\dobj}{\dGvobj}
\]
for the corresponding \Gvalg with underlying \Zsobj,
\[
 \isetelt \mapsto \Hom{\dcat}{\dobj}{\abs{\dGvobj}(\isetelt)}.
\]

The same conditions as in the ungraded case allow us to define an adjoint of the forgetful functor.

\begin{proposition}
 \label{prop:vfreegrade}
 Let \(\dcat\) be an (\(\m{E}\), \(\m{M}\)) category for some classes \(\m{E}\) of epimorphisms and \(\m{M}\) of monomorphisms.
 Suppose that:
\begin{enumerate}
\item \(\dcat\) has finite products,
\item \(\dcat\) is closed under filtered \co{}limits,
\item finite products commute with filtered \co{}limits,
\item \(\m{E}\) is closed under finite products, and
\item \(\dcat\) is \(\m{E}\)\hyp{}\co{well}\hyp{}powered.
\end{enumerate}
Then the forgetful functor \(\dGvcat \to \Zdcat\) has a left adjoint,
  \(\free{Gv} \colon \Zdcat \to \dGvcat\),
 which is called the \free{Gv}.
 \noproof
\end{proposition}

Exactly as in the ungraded case, we have \CGvalgobjs dual to \Gvalgobjs.

\begin{defn}
 Let \(\Gvcat\) be a variety of graded algebras, \(\dcat\) a category with finite \co{}products.
 A \emph{\dCGvobj{}} is a \OdGvobj.
 A \emph{morphism of \dCGvobjs} is a morphism in \(\Zdcat\) which intertwines the \dCGvobj structures.
 We denote \dCGvcat by \(\dCGvcat\).
\end{defn}

The analogue of Proposition~\ref{prop:conliftidgrade} is the following.

\begin{proposition}
 \label{prop:covliftgrade}
 To give an \Zdobj[\abs{\dCGvobj}] the structure of a \CGvalgobj is equivalent to giving a lift of the covariant hom\hyp{}functor
  \(\Hom{\dcat}{\abs{\dCGvobj}}{-} \colon \dcat \to \Zscat\)
 to a functor \(\dcat \to \Gvcat\).
 \noproof
\end{proposition}

Analogously to earlier,
 \(\Hom{\dcat}{\dCGvobj}{\dobj}\)
will denote the \Gvalg with underlying \Zsobj[\Hom{\dcat}{\abs{\dCGvobj}}{\dobj}].

\subsection{Adjunction Theorems}

The key tools for us are the links between representable functors and functors with adjoints.
The ungraded version of this result is a standard one from general algebra.

\begin{theorem}
 \label{th:adj}
 Let \(\dcat\) be a category that has finite products, is \co{complete}, is an (\(\m{E}\), \(\m{M}\)) category where \(\m{E}\) is closed under finite products, is \(\m{E}\)\hyp{}\co{well\hyp{}powered}, and is such that its finite products commute with filtered \co{limits}.
 Let \(\vcat\) be a variety of algebras.
 Let \(\fcat\) be a category with \co{equalisers}.
 Let
  \(\func{G} \colon \fcat \to \dvcat\)
 be a covariant functor.
 Then the following statements are equivalent.
 \begin{enumerate}
 \item
   \(\func{G}\) has a left adjoint.
 \item
  The composition
   \(\abs{\func{G}} \colon \fcat \to \dcat\)
  of \(\func{G}\) with the forgetful functor \(\dvcat \to \dcat\) has a left adjoint. \noproof
 \end{enumerate}
\end{theorem}

Working in the opposite category we obtain the following corollary.

\begin{corollary}
 \label{cor:adj}
 Let \(\dcat\) be a category that has finite products, is \co{complete}, is an (\(\m{E}\), \(\m{M}\)) category where \(\m{E}\) is closed under finite products, is \(\m{E}\)\hyp{}\co{well\hyp{}powered}, and is such that its finite products commute with filtered \co{limits}.
 Let \(\vcat\) be a variety of algebras.
 Let \(\fcat\) be a category with equalisers.
 Let
  \(\func{G} \colon \fcat \to \dvcat\)
 be a contravariant functor.
 Then the following statements are equivalent.
 \begin{enumerate}
 \item
   \(\func{G}\) is one of a mutually right adjoint pair.
 \item
  The composition
   \(\abs{\func{G}} \colon \fcat \to \dcat\)
  of \(\func{G}\) with the forgetful functor \(\dvcat \to \dcat\) is one of a mutually right adjoint pair.
  \noproof
 \end{enumerate}
\end{corollary}

Since a functor from a \co{complete} category into \(\scat\) is representable if and only if it has a left adjoint, using Proposition~\ref{prop:covlift} we deduce the following standard result of general algebra.

\begin{corollary}
 \label{cor:adjrep}
 \begin{enumerate}
 \item
  Let \(\dcat\) be a \co{complete} category, \(\vcat\) a variety of algebras.
  For a covariant functor
   \(\func{G} \colon \dcat \to \vcat\),
  the following statements are equivalent.
  \begin{enumerate}
  \item
    \(\func{G}\) has a left adjoint.
  \item
    \(\func{G}\) is representable by a \dCvobj.
  \item
    \(\abs{\func{G}}\) is representable by an \dobj.
  \end{enumerate}

 \item
  Let \(\dcat\) be a complete category, \(\vcat\) a variety of algebras.
  For a contravariant functor
   \(\func{G} \colon \dcat \to \vcat\),
  the following statements are equivalent.
  \begin{enumerate}
  \item
    \(\func{G}\) is one of a mutually right adjoint pair.
  \item
    \(\func{G}\) is representable by a \dvobj.
  \item
    \(\abs{\func{G}}\) is representable by an \dobj.
   \noproof
  \end{enumerate}
 \end{enumerate}
\end{corollary}

\medskip

These results generalise to the graded situtation without any substantial alterations.

\begin{theorem}
 \label{th:adjgrade}
 Let \(\dcat\) be a category that has finite products, is \co{complete}, is an (\(\m{E}\), \(\m{M}\)) category where \(\m{E}\) is closed under finite products, is \(\m{E}\)\hyp{}\co{well\hyp{}powered}, and is such that its finite products commute with filtered \co{limits}.
 Let \(\Gvcat\) be a variety of graded algebras.
 Let \(\fcat\) be a category with \co{equalisers}.
 Let
  \(\func{G} \colon \fcat \to \dGvcat\)
 be a covariant functor.
 Then the following statements are equivalent.
 \begin{enumerate}
 \item
   \(\func{G}\) has a left adjoint.
 \item
  The composition
   \(\abs{\func{G}} \colon \fcat \to \Zdcat\)
  of \(\func{G}\) with the forgetful functor \(\dGvcat \to \Zdcat\) has a left adjoint.
  \noproof
 \end{enumerate}
\end{theorem}

The graded version of Corollary~\ref{cor:adj} follows immediately.
To get the graded version of Corollary~\ref{cor:adjrep} we need to understand the relationship between graded adjunctions and representability.

\begin{lemma}
 \label{lem:gradjrep}
 Let \(\dcat\) be a \co{complete} category.
 A covariant functor
  \(\func{G} \colon \dcat \to \Zscat\)
 has a left adjoint if and only if it is representable by an \Zdobj.
\end{lemma}

\begin{proof}
 Suppose that \(\func{G}\) has a left adjoint, say
  \(\func{H} \colon \Zscat \to \dcat\).
 We extend this to a functor
  \(\func[mZo]{H} \colon \ZZscat \to \Zdcat\)
 in the obvious way.
 Let \(I\) be the \ZZsobj defined by
 \[
  \isetelt \mapsto \Bigg( \isetelt' \mapsto \begin{cases} \{*\} &
   \text{if } \isetelt = \isetelt' \\
   \emptyset &%
   \text{otherwise}
  \end{cases}
  \Bigg).
 \]
 Then for an \Zsobj[\Zsobj] we have 
 \begin{align*}
  \Hom{\Zscat}{I}{\Zsobj} &%
  = \big( \isetelt \mapsto \Hom{\Zscat}{I(\isetelt)}{\Zsobj} \big) \\
  &%
  \cong \big( \isetelt \mapsto \prod_{\isetelt' \in \iset} \Hom{\scat}{I(\isetelt)(\isetelt')}{\Zsobj(\isetelt')} \big) \\
  &%
  \cong \big(\isetelt \mapsto \Hom{\scat}{\{*\}}{\Zsobj(\isetelt)} \times \prod_{\isetelt' \ne \isetelt} \Hom{\scat}{\emptyset}{\Zsobj(\isetelt')} \big) \\
  &%
  \cong \big(\isetelt \mapsto \Zsobj(\isetelt) \times \prod_{\isetelt' \ne \isetelt} \{*\} \big) \\
  &%
  \cong \big(\isetelt \mapsto \Zsobj(\isetelt) \big) \\
  &%
  \cong \Zsobj,
 \end{align*}
 all natural in \(\Zsobj\).
 Hence there are natural isomorphisms
 \[
  \func{G}(\dobj) \cong \Hom{\Zscat}{I}{\func{G}(\dobj)} \cong \Hom{\dcat}{\func[mZo]{H}(I)}{\dobj}
 \]
 and so \(\func{G}\) is represented by the \Zdobj[\mZo{\func{H}}(I)].

 Conversely, suppose that \(\func{G}\) is represented by the \Zdobj[G].
 Let \(\Zsobj\) be an \Zsobj.
 We have the following natural isomorphisms of sets
 \begin{align*}
  \Hom{\Zscat}{\Zsobj}{\func{G}(\dobj)}
  &%
  \cong \Hom{\Zscat}{\Zsobj}{\Hom{\dcat}{G}{\dobj}} \\
  &%
  \cong \prod_{\isetelt \in \iset} \Hom{\scat}{\Zsobj(\isetelt)}{\Hom{\dcat}{G(\isetelt)}{\dobj}} \\
  &%
  \cong \prod_{\isetelt \in \iset} \Hom{\dcat}{G(\isetelt)}{\dobj}^{\Zsobj(\isetelt)} \\
  &%
  \cong \prod_{\isetelt \in \iset} \Hom{\dcat\big}{\coprod_{\Zsobj(\isetelt)} G(\isetelt)}{\dobj\big} \\
  &%
  \cong \Hom{\dcat\big}{\coprod_{\isetelt \in \iset}\coprod_{\Zsobj(\isetelt)} G(\isetelt)}{\dobj\big}.
 \end{align*}

 Therefore we define the functor
  \(\func{H} \colon \Zscat \to \dcat\)
 on objects by
 \[
  \func{H}(\Zsobj) = \coprod_{\isetelt \in \iset} \coprod_{\ZZsobj(\isetelt)} G(\isetelt)
 \]
 and in the obvious way on morphisms.
 This is the required left adjoint.
\end{proof}

As a corollary of this, and using Proposition~\ref{prop:conliftidgrade}, we deduce the graded version of Corollary~\ref{cor:adjrep}

\begin{corollary}
 \label{cor:adjrepgrade}
 \begin{enumerate}
 \item
  Let \(\dcat\) be a \co{complete} category, \(\Gvcat\) a variety of graded algebras.
  For a covariant functor
   \(\func{G} \colon \dcat \to \Gvcat\),
  the following statements are equivalent.
  \begin{enumerate}
  \item
    \(\func{G}\) has a left adjoint.
  \item
    \(\func{G}\) is representable by a \dCGvobj.
  \item
    \(\abs{\func{G}}\) is representable by an \Zdobj.
  \end{enumerate}

 \item
  Let \(\dcat\) be a complete category, \(\Gvcat\) a variety of graded algebras.
  For a contravariant functor
   \(\func{G} \colon \dcat \to \Gvcat\),
  the following statements are equivalent.
  \begin{enumerate}
  \item
    \(\func{G}\) is part of a mutually right adjoint pair.
  \item
    \(\func{G}\) is representable by a \dGvobj.
  \item
    \(\abs{\func{G}}\) is representable by an \Zdobj.
   \noproof
  \end{enumerate}
 \end{enumerate}
\end{corollary}

The work of the following sections can be viewed simply as applications of Corollaries~\ref{cor:adjrep} and~\ref{cor:adjrepgrade}.

\section{The Tall\enhyp{}Wraith Monoidal Structure}
\label{sec:twmon}

In this section we define and investigate a monoidal structure on \(\vCvcat\), and more generally on \(\GvCGvcat\).
As this monoidal structure was described for \(\vcat\) the category of commutative, unital rings in \cite{dtgw} we have elected to call it the \emph{Tall\enhyp{}Wraith monoidal structure}.

We give full details for the ungraded situation.
The graded analogues are recorded afterwards; they are straightforward modifications of the ungraded versions.

\begin{defn}
 Let \(\vcat\) be a variety of algebras.
 Let
  \(\covfun(\vcat, \vcat)\)
 denote the category of covariant functors from \(\vcat\) to itself.

 Let
  \(\covrep(\vcat, \vcat)\)
 be the full subcategory of
  \(\covfun(\vcat, \vcat)\)
 consisting of representable functors.
\end{defn}

By Proposition~\ref{prop:covlift}, a \vCvobj[\vCvobj] represents a covariant functor \(\vcat \to \vcat\) which we denote \(\cov{\vCvobj}\).
The standard Yoneda argument proves the following equivalence of categories.

\begin{lemma}
 The covariant functor
  \(\mOo{(\vCvcat)} \to \covrep(\vcat, \vcat)\),
 given on objects by
  \(\vCvobj \mapsto \cov{\vCvobj}\),
 is an equivalence of categories.
 \noproof
\end{lemma}

The category
 \(\covfun(\vcat, \vcat)\)
has an obvious monoidal structure coming from composition of functors.

\begin{proposition}
 \label{prop:covrepmon}
 Let \(\vcat\) be a variety of algebras.
 The category \(\covrep(\vcat,\vcat)\)
 inherits the monoidal structure of \(\covfun(\vcat,\vcat)\).
\end{proposition}

\begin{proof}
 To prove this we need to show that the composition of representable functors is again representable and that the identity functor is representable.

 As \(\vcat\) is \co{complete}, Corollary~\ref{cor:adjrep} applies with \(\dcat = \vcat\).
 Hence a covariant functor \(\vcat \to \vcat\) is representable if and only if it has a left adjoint.
 If two composable functors have left adjoints then their composition also has a left adjoint, whence the composition of two representable functors is again representable.

 Now let us consider the identity functor.
 We use Corollary~\ref{cor:adjrep} again to deduce that the identity functor on \(\vcat\) is representable if and only if the forgetful functor \(\vcat \to \scat\) is representable.
 The adjunction between the free and forgetful functors yields an isomorphism of sets, natural in \(\vobj\),
 \[
  \Hom{\vcat}{\free{v}(\{*\})}{\vobj} \cong \Hom{\scat}{\{*\}}{\sabs{\vobj}} \cong \sabs{\vobj}
 \]
 whence the forgetful functor is representable.
\end{proof}

From the equivalence of categories
 \(\vCvcat \cong \covrep(\vcat,\vcat)\)
we immediately deduce the following corollary.

\begin{corollary}
 The category \(\vCvcat\) has a monoidal structure such that the functor
  \(\mOo{(\vCvcat)} \to \covrep(\vcat, \vcat)\)
 is an equivalence of monoidal categories.
 \noproof
\end{corollary}

We shall write the monoidal structure as
 \((\vCvcat, \twprod, I)\).
It can be described explicitly as follows.

The objects in \vCvcat are objects in \vcat with extra structure.
This extra structure consists of \(\vcat\)\hyp{}morphisms from the underlying \vobj to iterated \co{products} of it.
Therefore any covariant functor
 \(\func{H} \colon \vcat \to \vcat\)
which preserves \co{products} has a natural lift to a functor
 \(\mCo{\func{H}} \colon \vCvcat \to \vCvcat\).
In particular, this is true if \(\func{H}\) is left adjoint to some functor
 \(\func{G} \colon \vcat \to \vcat\).
Hence for an \vobj[\vobj] and \vCvobj[\vCvobj] with underlying \vobj[\abs{\vCvobj}] the natural isomorphism of \sobjs
\[
 \Hom{\vcat}{\func{H}(\abs{\vCvobj})}{\vobj} \cong \Hom{\vcat}{\abs{\vCvobj}}{\func{G}(\vobj)}
\]
lifts naturally to a natural isomorphism of \vobjs
\[
 \Hom{\vcat}{\mCo{\func{H}}(\vCvobj)}{\vobj} \cong \Hom{\vcat}{\vCvobj}{\func{G}(\vobj)}.
\]

By Corollary~\ref{cor:adjrep}, the functor \(\cov{\vCvobj}\) has a left adjoint, which we denote by \(\ladj{\vCvobj}\), and this has a lift
 \(\Cladj{\vCvobj} \colon \vCvcat \to \vCvcat\)
as above.
For an \vobj[\vobj] and \vCvobjs[\vCvobj_1 \text{ and } \vCvobj_2] there are natural isomorphisms of \valgs
\begin{align*}
 \cov{(\Cladj{\vCvobj_2} (\vCvobj_1))}(\vobj) &%
 = \Hom{\vcat}{\Cladj{\vCvobj_2} (\vCvobj_1)}{\vobj} \\
 &%
 \cong \Hom{\vcat}{\vCvobj_1}{\cov{\vCvobj_2}(\vobj)} \\
 &%
 \cong \cov{\vCvobj_1}\cov{\vCvobj_2}(\vobj)
\end{align*}
and so
 \(\Cladj{\vCvobj_2}(\vCvobj_1)\)
represents the functor
 \(\cov{\vCvobj_1} \cov{\vCvobj_2}\).
By construction, the product
 \(\vCvobj_1 \twprod \vCvobj_2\)
also represents this functor and hence there is a natural isomorphism
\[
 \vCvobj_1 \twprod \vCvobj_2 \cong \Cladj{\vCvobj_2}(\vCvobj_1).
\]

The underlying \vobj of the unit, \(I\), of the monoidal structure is (isomorphic to) the free \valg on a singleton set.
The obvious isomorphisms of sets
\[
 \Hom{\vcat}{\free{v}(\{*\})}{\vobj} \cong \sabs{\vobj} \cong \sabs{\Hom{\vcat}{I}{\vobj}}
\]
underly isomorphisms of \valgs
\[
 \Hom{\vcat}{\free[hat]{v}(\{*\})}{\vobj} \cong \vobj \cong \Hom{\vcat}{I}{\vobj}
\]
and thus the isomorphism of \valgs,
 \(\sabs{I} \cong \free{v}(\{*\})\),
lifts to an isomorphism of \vCvobjs,
 \(I \cong \free[hat]{v}(\{*\})\).

\medskip

The above readily adapts to the following situations.

\begin{proposition}
 \label{prop:pairings}
 Let \(\dcat\) be a \co{complete} category, \(\vcat\) and \(\wcat\) varieties of algebras.
 There are products
 \begin{align*}
  \vCvcat \times \dCvcat &%
  \to \dCvcat, \\
  \mOo{(\vCvcat)} \times \dvcat &%
  \to \dvcat, \\
  \vCwcat \times \vCvcat &%
  \to \vCwcat,
 \end{align*}
 all compatible with the monoidal structure of \(\vCvcat\) and with composition of representable functors.
 We write all of the pairings using the notation \(-\twprod-\).
 \noproof
\end{proposition}

There are two things to note about this generalisation.
Firstly, special cases give two different pairings involving \(\vCvcat\) and \(\vcat\).
The first views \(\vcat\) as \(\svcat\) and so comes from the middle pairing above; in terms of functors we have
\[
 \con{(\vCvobj \twprod \vobj)}(\sobj) = \Hom{\scat}{\sobj}{\vCvobj \twprod \vobj} \cong \Hom{\vcat}{\vCvobj}{\Hom{\scat}{\sobj}{\vobj}}.
\]
The second views \(\vcat\) as \(\vCscat\) and so comes from the third pairing; in terms of functors we have
\[
 \cov{(\vobj \twprod \vCvobj)}(\vobj') = \Hom{\vcat}{\vobj \twprod \vCvobj}{\vobj'} \cong \Hom{\vcat}{\vobj}{\Hom{\vcat}{\vCvobj}{\vobj'}}.
\]
This latter pairing was the one considered in \cite{dtgw} in the case of commutative, unital rings.

The second thing to note about this generalisation is the annoyance of having a partially contravariant pairing.
Providing \(\dcat\) is sufficiently structured then this can be countered.

\begin{theorem}
 \label{th:shift}
 Let \(\dcat\) be a category satisfying the conditions of Theorem~\ref{th:adj}.
 Then there is a pairing
 \[
  \vCvcat \times \dvcat \to \dvcat, \qquad (\vCvobj, \dvobj) \mapsto \vCvobj \circledast \dvobj,
 \]
 which is covariant in both arguments and satisfies
 \[
  \Hom{\dvcat}{\vCvobj \circledast \dvobj}{\dvobj'} \cong \Hom{\dvcat}{\dvobj}{\vCvobj \twprod \dvobj'}
 \]
 naturally in all arguments.
\end{theorem}

\begin{proof}
 Let \(\vCvobj\) be a \vCvobj.
 Consider the covariant functor \(\dvcat \to \dvcat\) given by
  \(\dvobj \mapsto \vCvobj \twprod \dvobj\).
 We claim that this has a left adjoint.
 To show this, we shall prove that the functor \(\dvcat \to \dcat\),
  \(\dobj \mapsto \abs{\vCvobj \twprod \dvobj}\),
 has a left adjoint and then use Theorem~\ref{th:adj}.

 Firstly we observe that there is a pairing
  \(\vCvcat \times \dvCvcat \to \dvCvcat\)
 for which there is a natural isomorphism of \valgs
 \[
  \Hom{\vcat}{\vCvobj}{\Hom{\dvcat}{\dvCvobj}{\dvobj}} \cong \Hom{\dvcat}{\vCvobj \twprod \dvCvobj}{\dvobj}
 \]
 where \(\vCvobj\) is a \vCvobj, \(\dvCvobj\) is a \dvCvobj, and \(\dvobj\) is an \dvobjalt.

 Hence if \(\dobj\) is an \dobj, there is a natural isomorphism of \valgs
 \[
  \Hom{\dvcat}{\vCvobj \twprod \free[hat]{v}(\dobj)}{\dvobj} \cong \Hom{\vcat}{\vCvobj}{\Hom{\dvcat}{\free[hat]{v}(\dobj)}{\dvobj}} \cong \Hom{\vcat}{\vCvobj}{\Hom{\dcat}{\dobj}{\dvobj}} \cong \Hom{\dcat}{\dobj}{\vCvobj \twprod \dvobj}.
 \]
 In particular this produces a natural isomorphism on the underlying \sobjs.
 Its naturallity in all arguments demonstrates that the functor
  \(\dobj \mapsto \abs{\vCvobj \twprod \free[hat]{v}(\dobj)}\)
 is left adjoint to
  \(\dvobj \mapsto \abs{\vCvobj \twprod \dvobj}\).
 Thus by Theorem~\ref{th:adj}, the functor
  \(\dvobj \mapsto \vCvobj \twprod \dvobj\)
 has a left adjoint which we write
  \(\dvobj \mapsto \vCvobj \circledast \dvobj\).

 This is again natural in all arguments and so defines a pairing
 \[
  \vCvcat \times \dvcat \to \dvcat, \qquad (\vCvobj, \dvobj) \mapsto \vCvobj \circledast \dvobj.
 \]
 The adjunction proves the required identity.
\end{proof}

In a monoidal category it is natural to consider monoids.

\begin{defn}
 Let \(\vcat\) be a variety of algebras.
 A \emph{Tall\enhyp{}Wraith \(\vcat\)\hyp{}monoid} is a monoid in
  \((\vCvcat, \twprod, I)\).
 We write the category of such monoids as \(\vCvTtcat\).
\end{defn}

These were discussed briefly in \cite[\S 63, 64]{gbah}, though without explicit reference to the underlying monoidal structure on \(\vCvcat\).

Given a monoid one can consider modules for that monoid.
Since the monoidal category \(\vCvcat\) acts on other categories we can consider modules that are not \vCvobjs.
That is, if \(\vCvTtobj\) is a Tall\enhyp{}Wraith \(\vcat\)\hyp{}monoid and \(\dcat\) is a \co{complete} category then we can consider \dCvobjs[\dCvobj] for which there is a \(\dCvcat\)\hyp{}morphism
\[
 \vCvTtobj \twprod \dCvobj \to \dCvobj
\]
satisfying the required coherence conditions.

In \cite{gbah} the authors show that the category of \vobjs with an action of a Tall\enhyp{}Wraith monoid is again a variety of algebras.
Extending this, we easily see that a \dvobj or \dCvobj is a module for a Tall\enhyp{}Wraith \(\vcat\)\hyp{}monoid if and only if the corresponding functor \(\dcat \to \vcat\) factors through the category of \vobjs with an action of the Tall\enhyp{}Wraith \(\vcat\)\hyp{}monoid.

Two remarks are worth making at this juncture.
Firstly, if \(\wcat\) is another variety of algebras then the structure of a \(\vCvTtobj\)\hyp{}module on a \vCwobj does not have such an interpretation since a \vCwobj represents a functor \emph{out} of \vcat.
Secondly, due to the variance shift in the second pairing of Proposition~\ref{prop:pairings}, a \(\vCvTtobj\)\hyp{}module in \(\dvcat\) is better thought of as a \(\vCvTtobj\)\hyp{}\co{module} as the required morphism is
\[
 \dvobj \to \vCvTtobj \twprod \dvobj
\]
where the coherence condition involves the morphism \(\vCvTtobj \to \vCvTtobj \twprod \vCvTtobj\) in \(\mOo{\vCvcat}\)
We can surmount this using the product \(\circledast\) of Theorem~\ref{th:shift} since the adjunction turns a \co{action} as above into a more normal\hyp{}looking action.
That is to say, if \(\dvobj\) is a \(\vCvTtobj\)\hyp{}\co{}module for \(\twprod\) with action morphism
\[
 \dvobj \to \vCvTtobj \twprod \dvobj
\]
then it is a \(\vCvTtobj\)\hyp{}module for \(\circledast\) with action morphism
\[
 \vCvTtobj \circledast \dvobj \to \dvobj.
\]

\medskip

The adaptation of all this to the graded situation is straightforward.
The only part that may require a little extra explanation is the description of the monoidal structure.
A \GvCGvobj[\GvCGvobj] represents a functor
 \(\cov{\GvCGvobj} \colon \Gvcat \to \Gvcat\).
By Corollary~\ref{cor:adjrepgrade} this functor has a left adjoint
 \(\ladj{\GvCGvobj} \colon \Gvcat \to \Gvcat\)
which we extend to
 \(\Zladj{\GvCGvobj} \colon \ZGvcat \to \ZGvcat\)
in the obvious way.
This lift still preserves \co{products} and hence itself lifts to a functor
 \(\Cladj{\GvCGvobj} \colon \GvCGvcat \to \GvCGvcat\).
This has the property that the adjunction isomorphism lifts to an isomorphism of \Gvalgs
\[
 \Hom{\Gvcat}{\GvCGvobj_1}{\Hom{\Gvcat}{\GvCGvobj_2}{\Gvobj}} \cong \Hom{\Gvcat}{\Cladj{\GvCGvobj_2}(\GvCGvobj_1)}{\Gvobj}.
\]
Thus
 \(\Cladj{\GvCGvobj_2}(\GvCGvobj_1)\)
represents the functor
 \(\cov{\GvCGvobj_1} \cov{\GvCGvobj_2}\).

With this, we can generalise the results of this section to the graded sitation.  We record those generalisations here. 

\begin{proposition}
 Let \(\Gvcat\) be a variety of graded algebras.
 The category of representable functors,
  \(\covrep(\Gvcat, \Gvcat)\),
 inherits the monoidal structure of
  \(\covfun(\Gvcat, \Gvcat)\).
 \noproof
\end{proposition}

\begin{corollary}
 The category \(\GvCGvcat\) has a monoidal structure such that the functor
  \(\mOo{(\GvCGvcat)} \to \covrep(\Gvcat, \Gvcat)\)
 is an equivalence of monoidal categories.
 \noproof
\end{corollary}

We write this monoidal structure as
 \((\GvCGvcat, \twprod, I)\).

\begin{proposition}
 Let \(\dcat\) be a \co{complete} category, \(\Gvcat\) and \(\Gwcat\) varieties of graded algebras.
 There are products
 \begin{align*}
  \GvCGvcat \times \dCGvcat &%
  \to \dCGvcat, \\
  \mOo{(\GvCGvcat)} \times \dGvcat &%
  \to \dGvcat, \\
  \GvCGwcat \times \GvCGvcat &%
  \to \GvCGwcat,
 \end{align*}
 all compatible with the monoidal structure of \(\GvCGvcat\) and with composition of representable functors.
 We write all of the pairings using the notation \(-\twprod-\).
 \noproof
\end{proposition}

\begin{theorem}
 Let \(\dcat\) be a category satisfying the conditions of Theorem~\ref{th:adj}.
 Then there is a pairing
 \[
  \GvCGvcat \times \dGvcat \to \dGvcat, \qquad (\GvCGvobj, \dGvobj) \mapsto \GvCGvobj \circledast \dGvobj,
 \]
 which is covariant in both arguments and satisfies
 \[
  \Hom{\dGvcat}{\GvCGvobj \circledast \dGvobj}{\dGvobj'} \cong \Hom{\dGvcat}{\dGvobj}{\GvCGvobj \twprod \dGvobj'}
 \]
 naturally in all arguments.
 \noproof
\end{theorem}

\begin{defn}
 Let \(\Gvcat\) be a variety of graded algebras.
 A \emph{Tall\enhyp{}Wraith \(\Gvcat\)\hyp{}monoid} is a monoid in
  \((\GvCGvcat, \twprod, I)\).
 We write the category of such monoids as \(\GvCGvTtcat\).
\end{defn}

\medskip

\section{Free Algebras}
\label{sec:free}

In this section we shall prove Theorem~\ref{th:twmon}.
The first step is to show that the \free{v}
 \(\free{v} \colon \dcat \to \dvcat\)
lifts to a functor
 \(\free[hat]{v} \colon \dcat \to \dvCvcat\).
This is a straightforward application of Proposition~\ref{prop:covlift}.

\begin{theorem}
 \label{th:lift}
 Let \(\dcat\) be a \co{complete} category satisfying the conditions of Proposition~\ref{prop:vfree}.
 Let \(\vcat\) be a variety of algebras.
 There is a lift of
  \(\free{v} \colon \dcat \to \dvcat\)
 to a functor
  \(\free[hat]{v} \colon \dcat \to \dvCvcat\).
\end{theorem}

\begin{proof}
 Let \(\dobj\) be an \dobj and \(\dvobj\) a \dvobj.
 The natural isomorphism of sets
 \[
  \Hom{\dvcat}{\free{v}(\dobj)}{\dvobj} \cong \Hom{\dcat}{\dobj}{\abs{\dvobj}}
 \]
 shows that
  \(\Hom{\dvcat}{\free{v}(\dobj)}{\dvobj}\)
 is the underlying \sobj of a \valg.
 This structure is natural in both \(\dobj\) and \(\dvobj\) and hence defines a functor
  \(\func{\dobj} \colon \dvcat \to \vcat\).
 The composition of this functor with the forgetful functor, \(\vcat \to \scat\), is the functor
  \(\Hom{\dvcat}{\free{v}(\dobj)}{-}\)
 which is representable.
 Hence by Proposition~\ref{prop:covlift}, its representing object can be given the structure of a \vCvobj.
 The \valg structure on
  \(\Hom{\dvcat}{\free{v}(\dobj)}{\dvobj}\)
 is natural in \(\dobj\), whence the \Cvalg structure on \(\free{v}(\dobj)\) is functorial in \(\dobj\).
 This defines the required lift.
\end{proof}

Now let us consider the graded case.
The proof of the first part of Theorem~\ref{th:twmongrade} is similar to that of Theorem~\ref{th:twmon}.
To make the statement in its most general setting, we extend the definition of the diagonal functor.

\begin{defn}
 Let \(\dcat\) be a category with an initial object, \(\init{d}\).
 Define the \emph{diagonal} functor
  \(\diag \colon \dcat \to \ZZdcat\)
 on objects by
 \[
  \diag(\dobj)(\isetelt)(\isetelt') = \begin{cases} \dobj &%
   \text{if } \isetelt = \isetelt' \\
   \init{d} &%
   \text{otherwise},
  \end{cases}
 \]
 and similarly on morphisms.
\end{defn}

\begin{theorem}
 Let \(\dcat\) be a category satisfying the conditions of  Proposition~\ref{prop:vfreegrade}, \(\Gvcat\) a variety of graded algebras.
 The \free{Gv},
  \(\free{Gv} \colon \Zdcat \to \dGvcat\)
 extends to a functor
  \(\free[mZo]{Gv} \colon \ZZdcat \to \ZdGvcat\);
 the composition
  \(\free[mZo]{Gv} \diag \colon \dcat \to \ZdGvcat\)
 lifts to a functor
  \(\free[hat]{Gv} \colon \dcat \to \dGvCGvcat\).
\end{theorem}

\begin{proof}
 The first lift is the obvious one.
 We define
  \(\free[mZo]{Gv} \colon \ZZdcat \to \ZdGvcat\)
 on objects by
 \[
  \free[mZo]{Gv}\big( \isetelt \mapsto \Zdobj \big) = \big( \isetelt \mapsto \free{Gv}(\Zdobj) \big)
 \]
 and similarly on morphisms.

 For the second, let \(\dobj\) be an \dobj and \(\dGvobj\) a \dGvobj.
 Consider
 \[
  \Hom{\dGvcat}{\free[mZo]{Gv}\diag(\dobj)}{\dGvobj}.
 \]
 This is an \Zsobj.
 The \sobj corresponding to \(\isetelt \in \iset\) is
 \begin{align*}
  \Hom{\dGvcat}{\free[mZo]{Gv}\diag(\dobj)}{\dGvobj}(\isetelt) 
  &%
  = \Hom{\dGvcat}{(\free[mZo]{Gv}\diag(\dobj))(\isetelt)}{\dGvobj} \\
  &%
  = \Hom{\dGvcat}{\free{Gv}(\diag(\dobj)(\isetelt))}{\dGvobj} \\
  &%
  \cong \Hom{\Zdcat}{\diag(\dobj)(\isetelt)}{\abs{\dGvobj}} \\
  &%
  \cong \prod_{\isetelt' \in \iset} \Hom{\dcat}{\diag(\dobj)(\isetelt)(\isetelt')}{\abs{\dGvobj}(\isetelt')} \\
  &%
  \cong \Hom{\dcat}{\dobj}{\abs{\dGvobj}(\isetelt)}.
 \end{align*}
 The \Zsobj[\isetelt \mapsto \Hom{\dcat}{\dobj}{\abs{\dGvobj}(\isetelt)}] underlies a \Gvalg.
 The proof now proceeds as in that of Theorem~\ref{th:lift}.
\end{proof}

Notice that we do not get a lift of \(\free[mZo]{Gv}\) itself but only its composition with the diagonal functor.
There are special circumstances in which a lift of \(\free[mZo]{Gv}\) is possible, or is possible over a larger subcategory than the image of the diagonal functor, but in the general case the lift given above is the best possible.

\medskip

Let us return to the ungraded situation and complete the proof of Theorem~\ref{th:twmon} by examining the interaction between the \free{v} and the Tall\enhyp{}Wraith monoidal structure on \(\vCvcat\).

\begin{theorem}
 \label{th:mon}
 The lift of the \free{v},
  \(\free[hat]{v} \colon \scat \to \vCvcat\),
 is strong monoidal.
\end{theorem}

\begin{proof}
 Let \(\sobj_1\) and \(\sobj_2\) be \sobjs.
 We wish to show that there is a natural isomorphism
 \[
  \free[hat]{v}(\sobj_1 \times \sobj_2) \cong \free[hat]{v}(\sobj_1) \twprod \free[hat]{v}(\sobj_2),
 \]
 satisfying the standard coherence condition.

 Consider the functors \(\vcat \to \vcat\) represented by each side of the desired identity.
 Let \(\vobj\) be an \vobj; then we have natural isomorphisms of \valgs
 \begin{align*}
  \Hom{\vcat}{\free[hat]{v}(\sobj_1) \twprod \free[hat]{v}(\sobj_2)}{\vobj} 
  &%
  \cong \Hom{\vcat}{\free[hat]{v}(\sobj_1)}{\Hom{\vcat}{\free[hat]{v}(\sobj_2)}{\vobj}} \\
  &%
  \cong \Hom{\scat}{\sobj_1}{\Hom{\scat}{\sobj_2}{\vobj}} \\
  &%
  \cong \Hom{\scat}{\sobj_1 \times \sobj_2}{\vobj} \\
  &%
  \cong \Hom{\vcat}{\free[hat]{v}(\sobj_1 \times \sobj_2)}{\vobj}.
 \end{align*}
 Thus
  \(\free[hat]{v}(\sobj_1 \times \sobj_2)\)
 and
  \(\free[hat]{v}(\sobj_1) \twprod \free[hat]{v}(\sobj_2)\)
 represent the same functor \(\vcat \to \vcat\) and hence are isomorphic.
 This isomorphism is clearly natural in \(\sobj_1\) and \(\sobj_2\) and so defines the required natural isomorphism
 \[
  \free[hat]{v}(- \times - ) \cong \free[hat]{v}(-) \twprod \free[hat]{v}(-).
 \]
 The coherence condition is checked in the same way.

 We note that \(\{*\}\) is the unit for the monoidal structure on \(\scat\) and \(\free[hat]{v}(\{*\})\)
 is naturally isomorphic to the unit for the monoidal structure on \(\vCvcat\).
 Again, the coherence conditions are straightforward.
\end{proof}

Let us denote the category of monoids in
 \((\scat, \times, \{*\})\)
by \(\tcat\).

\begin{corollary}
 The functor
  \(\free[hat]{v} \colon \scat \to \vCvcat\)
 lifts to a functor
  \(\free[tilde]{v} \colon \tcat \to \vCvTtcat\).
 \noproof
\end{corollary}

\medskip

As before, we record the graded versions of these results.
The graded analogue of Theorem~\ref{th:mon} provides the rest of Theorem~\ref{th:twmongrade}.

\begin{proposition}
 The lift of the \free{Gv},
  \(\free[hat]{Gv} \colon \scat \to \GvCGvcat\),
 is strong monoidal.
 \noproof
\end{proposition}

\begin{corollary}
 The functor
  \(\free[hat]{Gv} \colon \scat \to \GvCGvcat\)
 lifts to a functor
  \(\free[tilde]{Gv} \colon \tcat \to \GvCGvTtcat\).
 \noproof
\end{corollary}

As we have been at pains to point out, the existence of the lifts of \(\free{v}\) and \(\free{Gv}\) are not hard to prove.
What is arguably difficult to see is where the \co{operation} morphisms come from on \(\free{v}(\dobj)\).
As we shall need this later in Theorem~\ref{th:cprop} and as it can be tricky to keep track of all the pieces, we shall expand on this now.
 
The \valg structure on
 \(\Hom{\dcat}{\dobj}{\dvobj}\)
comes entirely from the diagonal morphisms on \(\dobj\) and the \valg{}\hyp{}operations on \(\dvobj\), whereas that on
 \(\Hom{\dvcat}{\free[hat]{v}(\dobj)}{\dvobj}\)
must come from the \co{operations} on \(\free[hat]{v}(\dobj)\)
and the \co{diagonal} morphisms on \(\dvobj\).

Let us review the reason why
 \(\Hom{\dcat}{\dobj}{\dvobj}\)
is a \valg.
Let \(\oop\) be a \valg{}\hyp{}operation of arity \(n\).
Let
 \(\oop_{\dvobj} \colon \abs{\dvobj}^n \to \abs{\dvobj}\)
be the corresponding morphism in \(\dcat\).
Let
 \(f_1, \dotsc, f_n \colon \dobj \to \abs{\dvobj}\)
be \(\dcat\)\hyp{}morphisms.
The \(\dcat\)\hyp{}morphism
 \(\oop(f_1, \dotsc, f_n) \colon \dobj \to \abs{\dvobj}\)
is the composition
\[
 \dobj \xrightarrow{\Delta^n} \dobj^n \xrightarrow{f_1 \times \dotsb \times f_n} \abs{\dvobj}^n \xrightarrow{\oop_{\dvobj}} \abs{\dvobj}.
\]
Notice that all of the real information as to how to combine the morphisms \(f_1, \dotsc, f_n\) into a single morphism
 \(\oop(f_1, \dotsc, f_n)\)
is contained in the morphism \(\oop_{\dvobj}\) which, in the notation
 \(\Hom{\dcat}{\dobj}{\dvobj}\),
is encoded in the object on the \emph{right}.

The relevant diagram for a \dvCvobj[\dvCvobj] is that for \(\dvcat\)\hyp{}morphisms
 \(g_1, \dotsc, g_n \colon \abs{\dvCvobj} \to \dvobj\),
the \(\dvcat\)\hyp{}morphism
 \(\oop(g_1, \dotsc, g_n)\)
is the composition
\[
 \abs{\dvCvobj} \xrightarrow{\oop_{\dvCvobj}} \coprod_n \abs{\dvCvobj} \xrightarrow{g_1 \scoprod \dotsb \scoprod g_n} \coprod_n \dvobj \xrightarrow{{\mCo{\Delta}}^n} \dvobj.
\]
Thus in
 \(\Hom{\dvcat}{\dvCvobj}{\dvobj}\),
all of the real information as to how to combine the morphisms \(g_1, \dotsc, g_n\) into a single morphism
 \(\oop(g_1, \dotsc, g_n)\)
is encoded in the object on the \emph{left}.

What we shall show is that this view is subtly incorrect when the \dvCvobj is of the form \(\free[hat]{v}(\sobj)\).
What is really happening is that the \co{product} on \(\dvcat\) actually encodes all the \(\valg\)\hyp{}operations.
All that the \co{operation} on \(\free[hat]{v}(\sobj)\)
is doing is selecting the correct operation from the \co{product} by presenting the correct inputs.

To see how this works we need to examine the interaction between finite \co{product}s and products of \dvobjs.
We still have in mind an operation \(\oop\) of arity \(n\).
Let
 \(\dvobj_1, \dotsc, \dvobj_n\)
be \dvobjs.
Observe that there is a natural morphism in \(\dvcat\) given by the product of \co{}projections
\begin{equation}
 \label{eq:cotopro}
 \prod_{i=1}^n \dvobj_i \to \left(\coprod_{i=1}^n \dvobj_i \right)^n.
\end{equation}
Upon applying the forgetful functor we obtain a \(\dcat\)\hyp{}morphism which we compose with the following \(\dcat\)\hyp{}morphism from the operation \(\oop\)
\[
 \oop_{\scoprod \dvobj_i} \colon \abs{\left(\coprod_{i=1}^n \dvobj_i\right)^n} = \abs{\left(\coprod_{i=1}^n \dvobj_i\right)}^n \to \abs{\coprod_{i=1}^n \dvobj_i }
\]
to obtain a \(\dcat\)\hyp{}morphism
\[
 h_{\oop} \colon \prod_{i=1}^n \abs{\dvobj_i} = \abs{\prod_{i=1}^n \dvobj_i} \to \abs{\coprod_{i=1}^n \dvobj_i}.
\]
The label on the morphism is to remind us that we needed to have chosen an operation
 \(\oop \in \abs{\otype}\)
with arity \(n\) before we began this process and that \(h_\oop\) depends explicitly on this choice.

If we put all the \(\dvobj_i\) the same, say \(\dvobj\), then we can recover \(\oop_{\dvobj}\) from \(h_\oop\) by post\hyp{}composing with the underlying \(\dcat\)\hyp{}morphism of the \co{diagonal} in \(\dvcat\).
To see this, consider the diagram in Figure~\ref{fig:diag}.
\begin{figure}
\begin{centre}
\begin{tikzpicture}[column sep=30pt,row sep=30pt]
\matrix (diag) [matrix of math nodes]
{
\abs{\dvobj}^n \\
\abs{\left(\coprod_{i=1}^n \dvobj\right)^n} &
\abs{\dvobj}^n \\
\abs{\coprod_{i=1}^n \dvobj} &
\abs{\dvobj}. \\
};
\path[->] (diag-1-1) edge (diag-2-1)
          (diag-2-1) edge node[auto,swap] {\(\oop_{\scoprod \dvobj}\)} (diag-3-1)
          (diag-2-1) edge (diag-2-2)
          (diag-2-2) edge node[auto] {\(\oop_{\dvobj}\)} (diag-3-2)
          (diag-3-1) edge (diag-3-2)
          (diag-1-1.south west) [bend right=45] edge node[auto,swap] {\(h_\oop\)} (diag-3-1.north west);
\path[dashed,->] (diag-1-1.south east) edge (diag-2-2.north west);
\end{tikzpicture}
\end{centre}
\caption{Recovering \(\oop_{\dvobj}\) from \(h_\oop\).}
\label{fig:diag}
\end{figure}
The horizontal morphisms are induced by the \co{diagonal} \(\mCo{\Delta}\) in \(\dvcat\) and therefore the square commutes.
The composition on the left is \(h_\oop\) and on the right is \(\oop_{\dvobj}\).
Therefore we need to show that the composition
 \(\abs{\dvobj}^n \to \abs{(\coprod_{i=1}^n \dvobj)^n} \to \abs{\dvobj}^n\)
is the identity.
Let
 \(g_i \colon \dvobj \to \coprod_{i=1}^n \dvobj\)
be the \co{projection} onto the \(i\)th \co{factor}.
The morphism
 \(\dvobj^n \to (\coprod_{i=1}^n \dvobj)^n\)
is \(\prod_{i=1}^n g_i\).
The composition
 \(\mCo{\Delta} g_i \colon \dvobj \to \dvobj\)
is the identity, whence the composition
 \(\dvobj^n \to (\coprod_{i=1}^n \dvobj)^n \to \dvobj^n\)
is the product of the identity on each factor and hence is the identity morphism.

With \(h_\oop\) we can construct the \Cvalg structure on a free \dvobj.
Let \(\dobj\) be an \dobj.
There is a morphism in \(\dcat\), natural in \(\dobj\),
\[
 D \xrightarrow{\Delta} D^n \xrightarrow{\eta_{\dobj}^n} \abs{\free{v}(\dobj)}^n \xrightarrow{h_\oop} \abs{\coprod_{i=1}^n \free{v}(\dobj)},
\]
where
 \(\eta_{\dobj} \colon \dobj \to \abs{\free{v}(\dobj)}\)
is the unit of the adjunction.
Hence there is a \(\dvcat\)\hyp{}morphism
 \(\mCo{\oop}_{\dobj} \colon \free{v}(\dobj) \to \coprod_{i=1}^n \free{v}(\dobj)\)
such that
 \(\abs{\mCo{\oop}_{\dobj}} \eta_{\dobj}\)
is the above morphism.

Let \(\dvobj\) be a \dvobj.
Let
 \(f_i \colon \dobj \to \abs{\dvobj}\)
be morphisms in \(\dcat\), let
 \(\hat{f}_i \colon \free{v}(\dobj) \to \dvobj\)
be the corresponding morphisms of \dvobjs.
Consider the diagram in Figure~\ref{fig:coops} which takes place in \(\dcat\).
\begin{figure}
\begin{centre}
\begin{tikzpicture}[column sep=30pt,row sep=30pt]
\matrix (coop) [matrix of math nodes]
{
  \dobj &
   \prod_{i=1}^n \dobj & &
   \prod_{i=1}^n \abs{\dvobj} &
   \abs{\dvobj} \\
   \dobj &
   \prod_{i=1}^n &
   \prod_{i=1}^n \abs{\free{v}(\dobj)} &
   \prod_{i=1}^n \abs{\dvobj} &
   \abs{\dvobj} \\
   \abs{\free{v}(\dobj)} &
   &
   \abs{\coprod_{i=1}^n \free{v}(\dobj)} &
   \abs{\coprod_{i=1}^n \dvobj} &
   \abs{\dvobj} \\
};
\path[->,every node/.style={auto}]
 (coop-1-1) edge node {\(\Delta\)} (coop-1-2)
 (coop-1-1) edge node {\(=\)} (coop-2-1)
 (coop-1-2) edge node {\(\prod f_i\)} (coop-1-4)
 (coop-1-2) edge node {\(=\)} (coop-2-2)
 (coop-1-4) edge node {\(\oop_{\dvobj}\)} (coop-1-5)
 (coop-1-4) edge node {\(=\)} (coop-2-4)
 (coop-1-5) edge node {\(=\)} (coop-2-5)
 (coop-2-1) edge node {\(\Delta\)} (coop-2-2)
 (coop-2-1) edge node {\(\eta_{\dobj}\)} (coop-3-1)
 (coop-2-2) edge node {\(\prod \eta_{\dobj}\)} (coop-2-3)
 (coop-2-3) edge node {\(\prod \abs{\hat{f}_i}\)} (coop-2-4)
 (coop-2-3) edge node {\(h_\oop\)} (coop-3-3)
 (coop-2-4) edge node {\(\oop_{\dvobj}\)} (coop-2-5)
 (coop-2-4) edge node {\(h_\oop\)} (coop-3-4)
 (coop-2-5) edge node {\(=\)} (coop-3-5)
 (coop-3-1) edge node {\(\abs{\mCo{\oop}_{\dobj}}\)} (coop-3-3)
 (coop-3-3) edge node {\(\abs{\coprod \hat{f}_i}\)} (coop-3-4)
 (coop-3-4) edge node {\(\abs{\mCo{\Delta}}\)} (coop-3-5)
;
\end{tikzpicture}
\end{centre}
\caption{Co\hyp{}operations in free algebras}
\label{fig:coops}
\end{figure}
The outer two squares on the upper row obviously commute.
The middle square on the upper row commutes because
 \(f_i = \abs{\hat{f}_i}\eta_{\dobj}\)
for each \(i\).
The lower left square commutes by definition of \(\mCo{\oop}_{\dobj}\).
The lower middle square commutes by naturality of the morphism \(h_\oop\).
The square on the lower right commutes by the factorisation above of \(\oop_{\dvobj}\).
Hence the whole diagram commutes.

All of the objects and morphisms in the bottom row are the underlying objects and morphisms in \(\dcat\) of objects and morphisms in \(\dvcat\).
The corresponding diagram in \(\dvcat\)
\[
 \free{v}(\dobj) \xrightarrow{\mCo{\oop}_{\dobj}}
 \coprod_{i=1}^n \free{v}(\dobj) \xrightarrow{\coprod \hat{f}_i}
 \coprod_{i=1}^n \dvobj \xrightarrow{\mCo{\Delta}}
 \dvobj
\]
shows that the morphism \(\mCo{\oop}_{\dobj}\) is the required \co{operation} on \(\free{v}(\dobj)\).

Two things are important to note about the morphism \(\mCo{\oop}_{\dobj}\).
Firstly, it (or rather the morphism it is derived from) starts with the diagonal morphism \(\dobj \to \dobj^n\).
Secondly, the part that depends on \(\oop\) is the morphism
 \(\prod \abs{\free{v}(\dobj)} \to \abs{\coprod \free{v}(\dobj)}\);
this corresponds to ``selecting the right inputs''.

\section{Commutative Varieties of Algebras}
\label{sec:abelian}

In our motivating example, cohomology theories, the specific variety of algebras in question \emhyp{} graded commutative unital rings \emhyp{} is often viewed as being built on top of the category of abelian groups.
That is, rather than considering a ring as a ``set with structure'', it can be useful to regard it as an ``abelian group with structure'', where an abelian group is then a ``set with structure''.
This layering of varieties of algebras can be defined very generally.
An important issue is identifying characteristics of good midpoints.

One property of abelian groups that makes it a good choice is that it is a \emph{commutative} variety of algebras.
In an arbitrary variety of algebras the operations are defined on the underlying \sobj; in general they do not lift to the resulting category.
For abelian groups, though, this lift does occur.
That is, the operations defining abelian groups are themselves morphisms of abelian groups.

It is simple to write down the identities required for this to happen.
Each operation must commute with all other operations (including itself).
That is, if \(\vop_1\) and \(\vop_2\) are operations of arities \(n_1\) and \(n_2\) respectively then we require the identity
\begin{multline*}
  \vop_1(\vop_2(\velt_{1 1}, \dotsc, \velt_{1 n_2}), \dotsc, \vop_2(\velt_{n_1 1}, \dotsc, \velt_{n_1 n_2})) \\
=   \vop_2(\vop_1(\velt_{1 1}, \dotsc, \velt_{n_1 2}), \dotsc, \vop_1(\velt_{1 n_2}, \dotsc, \velt_{n_1 n_2}))
\end{multline*}

This is easier to see in ``matrix form''.
Let \(\valg\) be a variety of algebras and \(\vobj\) an \vobj.
Consider an \((n_1 \times n_2)\)\enhyp{}matrix of elements of \(\vobj\).
Using \(\vop_1\) and \(\vop_2\) there are two obvious ways of creating a new element of \(\vobj\).
The first way is to apply \(\vop_1\) to each column of the matrix resulting in a row vector of length \(n_2\) to which we apply \(\vop_2\).
The second way reverses this by applying \(\vop_2\) to each row of the matrix resulting in a column vector of length \(n_1\) to which we apply \(\vop_1\).
See Figure~\ref{fig:commops}.

\begin{figure}
\begin{centre}
\begin{tikzpicture}[column sep=25pt, row sep=30pt]
\matrix (commops) [matrix of math nodes]
{
\begin{bmatrix}
  \velt_{1 1} & \velt_{1 2} & \dotsb & \velt_{1 n_2} \\
  \velt_{2 1} & \velt_{2 2} & \dotsb & \velt_{2 n_2} \\
  \vdots & \vdots & \ddots & \vdots \\
  \velt_{n_1 1} & \velt_{n_1 2} & \dotsb & \velt_{n_1 n_2}
\end{bmatrix} &
\begin{bmatrix}
    \vop_2(\velt_{1 1}, \velt_{1 2}, \dotsc, \velt_{1 n_2}) \\
    \vop_2(\velt_{2 1}, \velt_{2 2}, \dotsc, \velt_{2 n_2}) \\
    \vdots \\
    \vop_2(\velt_{n_1 1}, \velt_{n_1 2}, \dotsc, \velt_{n_1 n_2})
\end{bmatrix} \\
\begin{bmatrix}
    \vop_1\begin{pmatrix}
    \velt_{1 1}\\ \velt_{2 1}\\ \vdots \\ \velt_{n_1 1}
    \end{pmatrix} &
    \vop_1\begin{pmatrix}
    \velt_{1 2} \\ \velt_{2 2} \\ \vdots \\ \velt_{n_1 2}
    \end{pmatrix} &
    \dotsb &
    \vop_1\begin{pmatrix}
    \velt_{1 n_2} \\ \velt_{2 n_2} \\ \vdots \\ \velt_{n_1 n_2}
    \end{pmatrix}
\end{bmatrix} \\
\vop_2
\left(
    \vop_1\begin{pmatrix}
    \velt_{1 1}\\ \velt_{2 1}\\ \vdots \\ \velt_{n_1 1}
    \end{pmatrix},
    \vop_1\begin{pmatrix}
    \velt_{1 2} \\ \velt_{2 2} \\ \vdots \\ \velt_{n_1 2}
    \end{pmatrix},
    \dotsc,
    \vop_1\begin{pmatrix}
    \velt_{1 n_2} \\ \velt_{2 n_2} \\ \vdots \\ \velt_{n_1 n_2}
    \end{pmatrix}
\right)
&\vop_1
\begin{pmatrix}
    \vop_2(\velt_{1 1}, \velt_{1 2}, \dotsc, \velt_{1 n_2}) \\
    \vop_2(\velt_{2 1}, \velt_{2 2}, \dotsc, \velt_{2 n_2}) \\
    \vdots \\
    \vop_2(\velt_{n_1 1}, \velt_{n_1 2}, \dotsc, \velt_{n_1 n_2})
\end{pmatrix} \\
 \\
};
\path[->,every node/.style={auto}]
 (commops-1-1) edge node {\(\vop_2\)} (commops-1-2)
 (commops-1-1) edge node {\(\vop_1\)} (commops-2-1)
 (commops-1-2) edge node {\(\vop_1\)} (commops-3-2)
 (commops-2-1) edge node {\(\vop_2\)} (commops-3-1)
;
\path
 (commops-3-1) -- node {\(=\)} (commops-3-2)
;
\end{tikzpicture}
\end{centre}
\caption{Commuting operations}
\label{fig:commops}
\end{figure}

Such varieties of algebras have turned up in the literature from time to time, the most relevant studies to our work being due to Freyd in \cite{pf2} and Linton \cite{fl2}.

\begin{defn}
 A \emph{commutative variety of algebras} is a variety of algebras in which every operation lifts to a morphism in the variety of algebras.
\end{defn}

It is interesting to note that the generalisation of this notion to graded algebras is not straightforward.
Operations in a graded variety of algebras are not morphisms of the underlying graded sets but are morphisms of components of such; therefore to define a commutative variety of graded algebras one would need to know that the operations came in suitable families which one could patch together to get a morphism of the underlying graded sets.
We shall not pursue this line here.

The category corresponding to a commutative variety of algebras has much structure which is explored in \cite{fl2}.
The following theorem summarises this structure.

\begin{theorem}[Linton]
Let \(\vcat\) be a commutative variety of algebras.
Then \(\vcat\) is a closed, symmetric monoidal category with identity (the) free \valg on a one element \sobj.

We shall write \((\vcat, \swprod, \free{v}(\{*\}))\) for the monoidal structure and \(\eHom{\vcat}{-}{-}\) for the enriched hom\hyp{}functor.

This structure has the property that the natural inclusion of sets
\[
  \Hom{\vcat}{-}{-} \to \Hom{\scat}{\abs{-}}{\abs{-}}
\]
lifts to a natural transformation of \(\vcat\)\enhyp{}valued functors
\[
  \eHom{\vcat}{-}{-} \to \Hom{\scat}{\abs{-}}{-}. \noproof
\]
\end{theorem}

The basic idea of this is quite simple.
If \(\vcat\) is a commutative variety of algebras then the \vops lift from \smors to \vmors and can be used to make each \valg a \vvobj and this defines the enriched hom\hyp{}functor.
The monodial product comes from extending the idea of bilinearity from vector spaces (or modules) to varieties of algebras, thus \(\Hom{\vcat}{\vobj_1 \swprod \vobj_2}{\vobj_3}\) is identified with the set of \smors \(\smor \colon \abs{\vobj_1} \times \abs{\vobj_2} \to \abs{\vobj_3}\) such that for each \(\velt \in \abs{\vobj_1}\), \(\smor(\velt,-) \colon \abs{\vobj_2} \to \abs{\vobj_3}\) underlies a \vmor, and similarly for each \(\velt \in \abs{\vobj_2}\).

This description of the structure focusses on the target in \(\Hom{\vcat}{\vobj_1}{\vobj_2}\).
It can be rephrased in terms of the source and the same argument that says that each \valg is naturally a \vvobj also says that each \valg is naturally a \vCvobj.

\begin{corollary}
 \label{cor:coab}
 Let \(\vcat\) be a commutative variety of algebras.
 There is a functor
  \(\func{C} \colon \vcat \to \vCvcat\)
 such that there is a natural isomorphism of \valgs
 \[
  \Hom{\vcat}{\func{C}(\vobj_1)}{\vobj_2} \cong \eHom{\vcat}{\vobj_1}{\vobj_2}.
 \]
\end{corollary}

\begin{proof}
 Let \(\vobj\) be an \vobj.
 Consider the functor
  \(\eHom{\vcat}{\vobj}{-} \colon \vcat \to \vcat\).
 The composition of this with the forgetful functor \(\vcat \to \scat\) is
  \(\Hom{\vcat}{\vobj}{-} \colon \vcat \to \scat\).
 This is represented by \(\vobj\) and so by Proposition~\ref{prop:covlift}, \(\vobj\) can be given the structure of a \vCvobj.
 As the assignment
  \(\vobj \to \eHom{\vcat}{\vobj}{-}\)
 is functorial, the \vCvobj structure of \(\vobj\) is also functorial in \(\vobj\).
 Hence we have a functor
  \(\func{C} \colon \vcat \to \vCvcat\)
 such that there is an isomorphism of \valgs
 \[
  \Hom{\vcat}{\func{C}(\vobj_1)}{\vobj_2} \cong \eHom{\vcat}{\vobj_1}{\vobj_2}
 \]
 naturally in all arguments.
\end{proof}

Some further properties of the functor \(\func{C}\) are given in the next theorem.

\begin{theorem}
 \label{th:cprop}
 The functor
  \(\func{C} \colon \vcat \to \vCvcat\)
 has the following properties.
 \begin{enumerate}
 \item
   \(\func{C}\) is right inverse to the forgetful functor \(\vCvcat \to \vcat\).

 \item
 \(\func{C}\) is a left adjoint (even, an enriched left adjoint).

 \item
  The lift of the \free{v},
   \(\free[hat]{v} \colon \scat \to \vCvcat\),
  factors as \(\func{C}\free{v}\).

 \item
   \(\func{C}\) is strong monoidal.
 \end{enumerate}
\end{theorem}

\begin{proof}
 \begin{enumerate}
 \item
  The functor \(\func{C}\) does not change the underlying \vobj and thus
   \(\abs{\func{C}} \colon \vcat \to \vcat\)
  is the identity functor.
  Hence \(\func{C}\) is a right inverse of the forgetful functor \(\vCvcat \to \vcat\).

  \item The right adjoint of \(\func{C}\) is the functor \(\func{V} \colon \vcat \to \vvcat\) which uses the lifts of the \vops to view \valgs as \vvobjs.
Since this is the implicit functor involved in the enriched hom\hyp{}functor, we have a natural isomorphism of \valgs
\[
  \eHom{\vcat}{\vobj_1}{\vobj_2} \cong \Hom{\vcat}{\vobj_1}{\func{V}(\vobj_2)}
\]
which, combined with the identity in \ref{cor:coab}, provides the (enriched) adjunction.

 \item
  For this part we need to show that the two ways of putting a \valg structure on
   \(\Hom{\vcat}{\free{v}(\sobj)}{\vobj}\)
  are the same.
  These two ways are
  \[
   \eHom{\vcat}{\free{v}(\sobj)}{\vobj} \text{ and } \Hom{\vcat}{\free[hat]{v}(\sobj)}{\vobj}.
  \]
  The \valg structure on the left hand side is characterised by the fact that the inclusion
  \[
   \Hom{\vcat}{\free{v}(\sobj)}{\vobj} \subseteq \Hom{\scat}{\abs{\free{v}(\sobj)}}{\abs{\vobj}}
  \]
  lifts to a morphism of \valgs
  \[
   \eHom{\vcat}{\free{v}(\sobj)}{\vobj} \subseteq \Hom{\scat}{\abs{\free{v}(\sobj)}}{\vobj}.
  \]

  That on the right is characterised by the fact that the isomorphism
  \[
   \Hom{\vcat}{\free{v}(\sobj)}{\vobj} \cong \Hom{\scat}{\sobj}{\abs{\vobj}}
  \]
  lifts to an isomorphism of \valgs
  \[
   \Hom{\vcat}{\free[hat]{v}(\sobj)}{\vobj} \cong \Hom{\scat}{\sobj}{\vobj}.
  \]

  Therefore we need to show that the morphism of \sobjs
  \[
   \Hom{\scat}{\sobj}{\abs{\vobj}} \cong \Hom{\vcat}{\free{v}(\sobj)}{\vobj} \subseteq \Hom{\scat}{\abs{\free{v}(\sobj)}}{\abs{\vobj}}
  \]
  lifts to a morphism of \valgs
  \[
   \Hom{\scat}{\sobj}{\vobj} \subseteq \Hom{\scat}{\abs{\free{v}(\sobj)}}{\vobj}.
  \]

  To show this we review the description of the \Cvalg structure on \(\free[hat]{v}(\sobj)\).
  Recall from Section~\ref{sec:free} that a \(\vcat\)\hyp{}operation \(\oop\) of arity \(n\) defines a \(\scat\)\hyp{}morphism 
  \[
   h_{\oop} \colon \prod_{i=1}^n \abs{\vobj_i} \to \abs{\coprod_{i=1}^n \vobj_i}.
  \]
  This is defined as the composition of two \smors.
  The first underlies a \(\vcat\)\hyp{}morphism and the second is a \(\vcat\)\hyp{}operation.
  Thus as, in our situation, \(\vcat\) is a commutative variety of algebras both of these components underlie \(\vcat\)\hyp{}morphisms and hence \(h_{\oop}\) underlies a \(\vcat\)\hyp{}morphism which we shall denote
  \[
   \hat{h}_{\oop} \colon \prod_{i=1}^n \vobj_i \to \coprod_{i=1}^n \vobj_i.
  \]

  Consider the lower left hand square of Figure~\ref{fig:coops}, reproduced in Figure~\ref{fig:lhcoop}.
\begin{figure}
\begin{centre}
\begin{tikzpicture}[column sep=30pt,row sep=30pt]
\matrix (lhcoop) [matrix of math nodes]
{
    \sobj &
    \prod_{i=1}^n \sobj &
    \prod_{i=1}^n \abs{\free{v}(\sobj)} \\
    \abs{\free{v}(\sobj)} &&
    \abs{\coprod_{i=1}^n \free{v}(\sobj)}\\
};
\path[->,every node/.style={auto}]
 (lhcoop-1-1) edge node {\(\Delta\)} (lhcoop-1-2)
 (lhcoop-1-1) edge node {\(\eta_{\sobj}\)} (lhcoop-2-1)
 (lhcoop-1-2) edge node {\(\eta_{\sobj}^n\)} (lhcoop-1-3)
 (lhcoop-1-3) edge node {\(h_{\oop}\)} (lhcoop-2-3)
 (lhcoop-2-1) edge node {\(\abs{\mCo{\oop}_{\sobj}}\)} (lhcoop-2-3)
;
\end{tikzpicture}
\end{centre}
\caption{Lower left-hand square of figure~\protect\ref{fig:coops}}
\label{fig:lhcoop}
\end{figure}
  The \(\vcat\)\hyp{}morphism \(\mCo{\oop}_{\sobj}\) is defined by the condition that this diagram commute.
  It is therefore given by
  \[
   \mCo{\oop}_{\sobj} \coloneqq \epsilon \free{v}(h_{\oop} \eta_{\sobj}^n \Delta)
  \]
  where \(\epsilon\) is the \co{unit} of the adjunction.
  We can simplify this expression using various facts.
  Firstly,
   \(\eta_{\sobj}^n \Delta = \Delta \eta_{\sobj}\).
  Secondly, \(h_{\oop}\) and \(\Delta\) underly \(\vcat\)\hyp{}morphisms and so we can use the formula
   \(\epsilon \free{v}(\abs{f}) = f \epsilon\).
  Finally,
   \(\epsilon \free{v}(\eta_{\sobj})\)
  is the identity on \(\free{v}(\sobj)\).
  Hence
  \[
   \mCo{\oop}_{\sobj} = \hat{h}_{\oop} \Delta.
  \]
  This means that with one slight modification we can lift the lower two lines of Figure~\ref{fig:coops} to a diagram in \(\vcat\).
\begin{figure}
\begin{centre}
\begin{tikzpicture}[column sep=30pt,row sep=30pt]
\matrix (vcoops) [matrix of math nodes]
{
    {\free{v}(\sobj)} &
    \prod_{i=1}^n {\free{v}(\sobj)} &
    \prod_{i=1}^n {\vobj} &
    {\vobj} \\
    {\free{v}(\sobj)} &
    {\coprod_{i=1}^n \free{v}(\sobj)} &
    {\coprod_{i=1}^n \vobj} &
    {\vobj} \\
 };
\path[->,every node/.style={auto}]
 (vcoops-1-1) edge node {\(\Delta\)} (vcoops-1-2)
 (vcoops-1-1) edge node {\(=\)} (vcoops-2-1)
 (vcoops-1-2) edge node {\(\prod {\hat{f}_i}\)} (vcoops-1-3)
 (vcoops-1-2) edge node {\(\hat{h}_\oop\)} (vcoops-2-2)
 (vcoops-1-3) edge node {\(\oop_{\vobj}\)} (vcoops-1-4)
 (vcoops-1-3) edge node {\(\hat{h}_\oop\)} (vcoops-2-3)
 (vcoops-1-4) edge node {\(=\)} (vcoops-2-4)
 (vcoops-2-1) edge node {\(\mCo{\oop}_{\sobj}\)} (vcoops-2-2)
 (vcoops-2-2) edge node {\(\coprod \hat{f}_i\)} (vcoops-2-3)
 (vcoops-2-3) edge node {\(\mCo{\Delta}\)} (vcoops-2-4)
;
\end{tikzpicture}
\end{centre}
\caption{Lower lines of figure~\protect\ref{fig:coops} lifted to \(\vcat\)}
\end{figure}
  The lower line gives the \valg structure on
   \(\Hom{\vcat}{\free{v}(\sobj)}{\vobj}\)
  defined by the lift \(\free[hat]{v}(\sobj)\)
  and the upper line gives the \valg structure coming from the enrichment.
  The two therefore agree and hence
   \(\free[hat]{v} = \func{C}\free{v}\).

 \item
  As
   \((\vcat, \swprod, \free{v}(\{*\})\)
  is a closed symmetric monoidal category, for \vobjs[\vobj_1, \vobj_2, \vobj_3] there is an isomorphism of \valgs natural in all arguments
  \[
   \eHom{\vcat}{\vobj_1}{\eHom{\vcat}{\vobj_2}{\vobj_3}} \cong \eHom{\vcat}{\vobj_1 \swprod \vobj_2}{\vobj_3}.
  \]
  Hence there are natural isomorphisms of \valgs
  \begin{align*}
   \Hom{\vcat}{\func{C}(\vobj_1 \swprod \vobj_2)}{\vobj_3} 
   &%
   \cong \eHom{\vcat}{\vobj_1 \swprod \vobj_2}{\vobj_3} \\
   &%
   \cong \eHom{\vcat}{\vobj_1}{\eHom{\vcat}{\vobj_2}{\vobj_3}} \\
   &%
   \cong \Hom{\vcat}{\func{C}(\vobj_1)}{\Hom{\vcat}{\func{C}(\vobj_2)}{\vobj_3}} \\
   &%
   \cong \Hom{\vcat}{\func{C}(\vobj_1) \twprod \func{C}(\vobj_2)}{\vobj_3}.
  \end{align*}
  Thus
   \(\func{C}(\vobj_1 \swprod \vobj_2)\)
  and
   \(\func{C}(\vobj_1) \twprod \func{C}(\vobj_2)\)
  represent the same functor \(\vcat \to \vcat\) whence are isomorphic as \vCvobjs.
  Moreover, this isomorphism is natural in both \(\vobj_1\) and \(\vobj_2\).
  The associativity coherence follows by a similar argument.
  From above we have an isomorphism
   \(\func{C}\free{v}(\{*\}) \cong \free[hat]{v}(\{*\})\)
  and it is easy to show that this also satisfies the required coherences.
  \qedhere
 \end{enumerate}
\end{proof}

As an immediate corollary we obtain the result that \(\free{v}\) itself is strong monoidal.

\begin{corollary}
 Let \(\vcat\) be a commutative variety of algebras.
 The functor
  \(\free{v} \colon \scat \to \vcat\)
 is strong symmetric monoidal.
\end{corollary}

\begin{proof}
 Let \(\sobj_1\) and \(\sobj_2\) be \sobjs.
 Then
 \begin{align*}
  \free{v}(\sobj_1 \times \sobj_2) &%
  = \abs{\free[hat]{v}(\sobj_1 \times \sobj_2)} \\
  &%
  \cong \abs{\free[hat]{v}(\sobj_1) \twprod \free[hat]{v}(\sobj_2)} \\
  &%
  \cong \abs{\func{C}\free{v}(\sobj_1) \twprod \func{C}\free{v}(\sobj_2)} \\
  &%
  \cong \abs{\func{C}(\free{v}(\sobj_1) \swprod \free{v}(\sobj_2))} \\
  &%
  \cong \free{v}(\sobj_1) \swprod \free{v}(\sobj_2).
 \end{align*}
 Also, \(\free{v}(\{*\})\) is the unit for the monoidal structure on \(\vcat\).
 The coherences follow by similar arguments.
\end{proof}

\begin{examples}
 \begin{enumerate}
 \item
  Let \(\vcat\) be the category of abelian groups.
  It is not hard to check that an abelian group \(\vobj\) admits precisely one \co{}\vcat structure given by:
  \begin{align*}
   \mCo{\alpha} &%
   \colon \vobj \to \vobj \scoprod \vobj, &
   \velt &%
   \mapsto (\velt, \velt), \\
   \mCo{\nu} &%
   \colon \vobj \to \vobj, & \velt &%
   \mapsto - \velt, \\
   \mCo{\zeta} &%
   \colon \vobj \to \{0\}, & \velt &%
   \mapsto 0.
  \end{align*}
  Hence \(\vCvcat \cong \vcat\).
  As the functor \(\vcat \to \vCvcat\) is \(\func{C}\), this is an isomorphism of symmetric monoidal categories.

 \item
  Let \(R\) be a commutative, unital ring and \(\vcat\) the category of \(R\)\hyp{}modules.
  The functor \(\func{C}\) assigns to an \(R\)\hyp{}module the obvious \(R\)\hyp{}bimodule.
  This need not be an equivalence of categories since an \(R\)\hyp{}module may admit several distinct structures as an \(R\)\hyp{}bimodule.
  For example, let \(R = \C\) and define a right \(\C\)\hyp{}action on a left \(\C\)\hyp{}module \(\vobj\) by
   \(\velt \cdot \lambda = \overline{\lambda} \velt\).

  This example shows that \(\vCvcat\) need not be \emph{symmetric} monoidal, even if \(\vcat\) is a commutative variety of algebras, and that the forgetful functor \(\vCvcat \to \vcat\) need not be monoidal.
 \end{enumerate}
\end{examples}

\medskip

Let us write \(\vStcat\) for the category of monoids in \(\vcat\).
Recall that we write \(\tcat\) for the category of monoids in \(\scat\).

\begin{corollary}
 Let \(\vcat\) be a commutative variety of algebras.
 The \free{v},
  \(\free{v} \colon \scat \to \vcat\)
 lifts to a functor
  \(\free[check]{v} \colon \tcat \to \vStcat\).
 \noproof
\end{corollary}

\medskip

One can do general algebra in a symmetric monoidal category by using the monoidal product instead of the cartesian product and by using PROPs instead of algebraic theories.
Thus  we can consider the category of algebras of a given PROP in \((\vcat, \swprod, \free{v}(\{*\}))\) for \(\vcat\) a commutative variety of algebras.
When we specialise to a \scat\/\enhyp{}based operad, we obtain the following result.

\begin{theorem}
 Let \(\vcat\) be a commutative variety of algebras, let \(\wcat\) be a \scat\/\enhyp{}based operad.
 Let \(\vSwcat\) be \vSwcat.
 Then \(\vSwcat\) is a variety of algebras built by adding operations and identities to (a presentation of) \(\vcat\).
\end{theorem}

\begin{proof}
The new operations are morphisms \(\wmor \colon \vobj \swprod \dotsb \swprod \vobj \to \vobj\).
The characterisation of \(\vobj \swprod \dotsb \swprod \vobj\) as representing \(k\)\enhyp{}fold \(\vcat\)\enhyp{}homomorphisms with source \(\abs{\vobj} \times \dotsb \times \abs{\vobj}\) means that \(\wmor\) defines a \smor
\[
  \abs{\wmor} \colon \abs{\vobj} \times \dotsb \times \abs{\vobj} \to \abs{\vobj}.
\]
Thus we have a new operation.
This argument repeats to show that the identities satisfied by the operations in the operad \wcat are also satisfied by the new operations.

All that is left is to encode the fact that the new operations came from an operad.
This means that a new operation, say \(\abs{\wmor}\), lifts to a \vmor in each variable in turn.
We can encode that in identities as
\begin{multline*}
  \abs{\wmor}\big(\vop(\velt_1, \dotsc, \velt_n), \velt'_2, \dotsc, \velt'_m\big) \\
 = \vop\big(\abs{\wmor}(\velt_1, \velt'_2, \dotsc, \velt'_m),\abs{\wmor}(\velt_2, \velt'_2, \dotsc, \velt'_m),\dotsc,\abs{\wmor}(\velt_n, \velt'_2, \dotsc, \velt'_m)\big)
\end{multline*}
and similarly in the other entries.

When added to the operations and identities of \(\vcat\), these operations and identities make \(\vSwcat\) into a variety of algebras.
\end{proof}

Although we think of the relationship between \(\vSwcat\) and \(\vcat\) as similar to that between, say, \(\vcat\) and \(\scat\), there are important differences.
The main one is that whilst
 \(\Hom{\scat}{\sobj}{\abs{\vobj}}\)
underlies a \valg for any \sobj[\sobj]
 \(\Hom{\vcat}{\vobj}{\abs{\vSwobj}}\)
is not necessarily a \walg.
What may be missing is an analogue of the diagonal morphisms such as
 \(\sobj \to \sobj \times \sobj\).

\begin{defn}
 Let \(\vcat\) be a commutative variety of algebras.
 A \vCStobj is an \vobj[\vCStobj] together with \(\vcat\)\hyp{}morphisms
  \(\Delta \colon \vCStobj \to \vCStobj \swprod \vCStobj\)
 and
  \(\epsilon \colon \vCStobj \to \free{v}(\{*\})\)
 such that \(\Delta\) is \co{associative} and \(\epsilon\) is a \co{unit} for \(\Delta\).

 We denote \vCStcat by \(\vCStcat\).
\end{defn}

\begin{lemma}
The functor \(\vCStcat \times \vSwcat \to \scat\) defined by
\[
  (\vCStobj, \vSwobj) \mapsto \Hom{\vcat}{\vCStobj}{\vSwobj}
\]
lifts to a functor \(\vCStcat \times \vSwcat \to \wcat\).
\noproof
\end{lemma}

\vCStcatu inherits a symmetric monoidal structure from \(\vcat\).

As \(\vSwcat\) is a variety of algebras in the usual sense, the functor \(\vSwcat \to \vcat\) preserves underlying \sobjs and thus has a left adjoint,
 \(\free{vSwv} \colon \vcat \to \vSwcat\),
which interacts with the free functors from \(\scat\) in that
 \(\free{vSw} \cong \free{vSwv}\free{v}\).

The following is a straightforward adaptation of Theorem~\ref{th:twmon}.

\begin{proposition}
 Let \(\vcat\) be a commutative variety of algebras, \(\wcat\) an operad.
 Let \(\ucat = \vSwcat\).
 The \free{u} from \(\vcat\),
  \(\free{uv} \colon \vcat \to \ucat\),
 lifts to a functor
  \(\free[hat]{uv} \colon \vCStcat \to \uCucat\).
 This functor is strong monoidal.

 The functor
  \(\free{v} \colon \scat \to \vcat\)
 factors through \(\vCStcat\) and the functor
  \(\free[hat]{u} \colon \scat \to \uCucat\)
 factors as
  \(\free[hat]{u} \cong \free[hat]{uv} \free{v}\).
 \noproof
\end{proposition}

For the purpose of finding Tall\enhyp{}Wraith monoids, the fact that \(\free[hat]{u} \colon \vCStcat \to \uCucat\) is strong monoidal means that it lifts to a functor \(\vCStStcat \to \uCuStcat\), where \(\vCStStcat\) is \vCStStcat and \(\uCuStcat\) is \uCuStcat (that is, the category of Tall\enhyp{}Wraith \(\ucat\)\enhyp{}monoids). 

\begin{example}
 Let \(\vcat\) be the category of abelian groups.
 The category of commutative, unital rings is of the form \(\vSwcat\) with \(\wcat\) the operad corresponding to commutative monoids.
 Thus the free commutative ring functor from abelian groups to rings defines functors from the category of \co{algebras} (in the traditional sense) to the category of birings and from the category of (not necessarily commutative) bialgebras to the category of biring triples.
 These are specific instances of the examples given in \cite{jbbw}.
\end{example}

Figure~\ref{fig:lifts} summarises this structure.

\begin{figure}
\begin{centre}
\begin{tikzpicture}
[column sep=30pt,row sep=30pt]
\matrix (free) [matrix of math nodes]
{
&&& \uCuStcat \\
&&& \uCucat \\
&&& \ucat \\
\vCStStcat & \vCStcat & \vcat \\
&&&& \scat \\
&&&&& \tcat \\
};
\path[->,every node/.style={auto}]
 (free-1-4) edge (free-2-4)
 (free-2-4) edge (free-3-4)
 (free-4-1) edge (free-4-2)
 (free-4-2) edge (free-4-3)
 (free-5-5) edge node[swap,pos=.7] {\(\free{v}\)} (free-4-3)
 (free-5-5) edge node {\(\free[hat]{v}\)} (free-4-2)
 (free-6-6) edge node {\(\free[tilde]{v}\)} (free-4-1)
 (free-5-5) edge node[pos=.7] {\(\free{u}\)} (free-3-4)
 (free-5-5) edge node[swap] {\(\free[hat]{u}\)} (free-2-4)
 (free-6-6) edge node[swap] {\(\free[tilde]{u}\)} (free-1-4)
 (free-6-6) edge (free-5-5)
 (free-4-1) edge node {\(\free[tilde]{uv}\)} (free-1-4)
 (free-4-2) edge node {\(\free[hat]{uv}\)} (free-2-4)
 (free-4-3) edge node {\(\free{uv}\)} (free-3-4)
;
\end{tikzpicture}
\end{centre}
\caption{Lifts of free functors}
\label{fig:lifts}
\end{figure}

\medskip

We said earlier that there was not a simple definition of a commutative variety of \emph{graded} algebras.
On the other hand, the immediately preceding results have used two theories, a commutative variety of algebras and an operad.
It is entirely possible to replace the operad by a graded one (which is more commonly known as a coloured operad).
Exactly the same results hold.

\begin{proposition}
 Let \(\vcat\) be a commutative variety of algebras, let \(\Gwcat\) be a graded operad.
 Let \(\vSGwcat\) be \vSGwcat.
 Then \(\vSGwcat\) is a variety of graded algebras.
 \noproof
\end{proposition}

\begin{proposition}
 Let \(\vcat\) be a commutative variety of algebras, \(\Gwcat\) a graded operad.
 Let \(\Gucat = \vSGwcat\).
 The composition of
  \(\free[mZo]{Guv} \colon \ZZvcat \to \ZGucat\)
 with the diagonal functor
  \(\diag \colon \vcat \to \ZZvcat\)
 lifts to a functor
  \(\free[hat]{Guv} \colon \vCStcat \to \GuCGucat\).
 This functor is strong monoidal.

 The functor
  \(\free[hat]{Gu} \colon \scat \to \GuCGucat\)
 factors as
  \(\free[hat]{Gu} \cong \free[hat]{Guv} \free{v}\).
 \noproof
\end{proposition}

\section{Toy Cohomology Theories}
\label{sec:toys}

As mentioned in the introduction, our interest in these structures comes from studying unstable operations of cohomology theories.
The situation for cohomology theories is somewhat complicated both by the grading and the filtrations.
In this section we shall consider some ``toy cohomology theories'' in which, by simplifying wildly, we can remove both the grading and filtration and reveal the simplicity of the theory.
Of course, these are not true cohomology theories, but they retain the basic structure needed to illustrate our purpose in introducing Tall\enhyp{}Wraith monoids to describe the structure of cohomology operations.

Removing the grading is easy: we work with a ring object rather than a graded ring object.
To remove the filtration, we drop the topology.
Thus we consider ring objects in \(\scat\); in other words, rings.
Actually, we consider commutative, unital rings since cohomology theories are (graded) commutative and unital.
Henceforth, a ``ring'' will be assumed to be commutative and unital,  and we shall write \(\rcat\) for \rcat.

In \cite{MR2559638}, the cohomology theories that we studied satisfied the K\"unneth formula.
In our simplified situation, we need the map
\[
 \otimes_{\robj}^n \Hom{\scat}{\abs{\robj}}{\robj} \to  \Hom{\scat}{\abs{\robj}^n}{\robj}
\]
to be an isomorphism.
Here, \(\robj\) is a ring and \(n \in \N\).

This map is constructed as follows.
Each projection map \(\abs{\robj}^n \to \abs{\robj}\) defines an \rmor
\[
  \Hom{\scat}{\abs{\robj}}{\robj} \to \Hom{\scat}{\abs{\robj}^n}{\robj}.
\]
Putting these together, we obtain an \rmor
\[
 \coprod_n \Hom{\scat}{\abs{\robj}}{\robj} \to \Hom{\scat}{\abs{\robj}^n}{\robj}.
\]

At first sight, the ring \(\robj\) defines a functor \(\con{\robj} \coloneqq \Hom{\scat}{-}{\robj} \colon \scat \to \rcat\).
However, the rings in the image of this functor are actually \bobjs.
Instead, therefore, we work in \bcat, say \(\bcat\), and the correct coproduct is the tensor product over \(\robj\).
Thus we obtain an \bmor:
\begin{equation}
\label{eq:copiso}
 \otimes_{\robj}^n \Hom{\scat}{\abs{\robj}}{\robj} \to  \Hom{\scat}{\abs{\robj}^n}{\robj}.
\end{equation}

The reason that we want this to be an isomorphism is because we want to define a \co{ring} structure on \(\Hom{\scat}{\abs{\robj}}{\robj}\).
That is, for each \bop, say \(\bop\), of arity \(n\), we need a \bmor
\[
  \Hom{\scat}{\abs{\robj}}{\robj} \to \coprod_n \Hom{\scat}{\abs{\robj}}{\robj}.
\]
The \bop is a function \(\abs{\robj}^n \to \abs{\robj}\) and so defines a \bmor
\[
  \Hom{\scat}{\abs{\robj}}{\robj} \to \Hom{\scat}{\abs{\robj}^n}{\robj}.
\]
If the \rmors from \eqref{eq:copiso} are isomorphisms, we obtain the required morphism
\[
  \Hom{\scat}{\abs{\robj}}{\robj} \to \Hom{\scat}{\abs{\robj}^n}{\robj} \cong  \otimes_{\robj}^n \Hom{\scat}{\abs{\robj}}{\robj} \cong \coprod_n \Hom{\scat}{\abs{\robj}}{\robj}.
\]

To get isomorphisms in \eqref{eq:copiso} for \(n \in \N\), we need \(\robj\) to be finite.

\begin{defn}
By a \emph{toy cohomology theory} we mean the following situation.
We have a finite ring \(\robj\) which we use to define a contravariant functor from \scat to \bcat.
\end{defn}

The analogues of the main results of \cite{MR2559638} for this toy example are recorded in the following theorem.

\begin{theorem}
The monoidal structure on \(\Hom{\scat}{\abs{\robj}}{\abs{\robj}}\) makes \(\Hom{\scat}{\robj}{\robj}\) into a Tall\enhyp{}Wraith \(\bcat\)\enhyp{}monoid.

The functor \(\con{\robj} \colon \scat \to \bcat\) lifts to a functor to the category of \(\Hom{\scat}{\robj}{\robj}\)\enhyp{}modules in \bcat.
\end{theorem}

As in \cite{MR2559638}, we could also examine operations between two different toy cohomology theories.
As this is intended purely as toy examples, we concentrate on the main ideas and leave this to another time.

\section{Examples of Tall\enhyp{}Wraith Monoids}
\label{sec:examples}

Our primary purpose in writing this paper is to exhibit examples of Tall\enhyp{}Wraith monoids.
We have three routes available to us.
For some varieties of algebras, in particular (but not limited to) commutative varieties of algebras, the category of Tall\enhyp{}Wraith monoids is known and the corresponding Tall\enhyp{}Wraith monoids are similarly known.
In this case, we simply identify the categories.
Secondly, we can construct Tall\enhyp{}Wraith monoids using free functors.
This relies on having a simpler variety of algebras where the Tall\enhyp{}Wraith monoidal structure is known and wherein we consider Tall\enhyp{}Wraith monoids as simple objects (at least, simple relative to the variety of algebras that we wish to study).
Thirdly, we can construct Tall\enhyp{}Wraith monoids from toy cohomology theories.
In this section we shall follow those routes and give examples of Tall\enhyp{}Wraith monoids.

\subsection{Known Tall\enhyp{}Wraith Monoids}

For some varieties of algebras, the Tall\enhyp{}Wraith monoidal structure is already known, albeit not by that name.

\begin{enumerate}
\item
  Let \(\vcat = \scat\).
  Since \(\vcat\) has a presentation with no operations or identities, it is clear that \(\vCvcat = \scat\) and the Tall\enhyp{}Wraith monoidal structure is simply
   \((\scat, \times, \{*\})\).
Thus Tall\enhyp{}Wraith monoids in \(\scat\) are simply monoids.

 \item
  Let \(\vcat\) be the category of abelian groups.
  We saw in Section~\ref{sec:abelian} that \(\vCvcat\) is isomorphic to \(\vcat\) and under this isomorphism the Tall\enhyp{}Wraith monoidal structure is simply \((\vcat, \otimes, \Z)\).
Thus a Tall\enhyp{}Wraith monoid in \(\vcat\) is a unital ring.

 \item
  Let \(\vcat\) be the category of groups.
  It is an old theorem of D.\ M.\ Kan, \cite{MR0111035}, that the underlying group of a \co{group} in \(\vcat\) is a free group and that the functor
   \(\free[hat]{v} \colon \scat \to \vCvcat\)
  is an equivalence of monoidal categories.
A Tall\enhyp{}Wraith monoid in groups consists of a free group, a choice of basis for that free group, and a monoid structure on that basis.

 \item
  Let \(R\) be a unital ring and let \(\vcat\) be the category of left \(R\)\hyp{}modules.
  Then \(\vCvcat\) is the category of \(R\)\hyp{}bimodules and the monoidal structure is
   \((\vCvcat, {}_R \otimes_R, R)\).
Thus a Tall\enhyp{}Wraith \(\vcat\)\hyp{}monoid is an \(R\)\hyp{}algebra: the compatibility conditions for the unit and multiplication ensure that \(R\) acts centrally.

As this is such a well\hyp{}studied case, it is worth taking a moment to identify some of the other pieces of our structure with their classical counterparts. 

  Let \(R'\) be another unital ring and let \(\wcat\) be the category of left \(R'\)\hyp{}modules.
  Then \(\wCvcat\) is the category of \(R'\)\hyp{}\(R\)\hyp{}bimodules.
  The first pairing of Proposition~\ref{prop:pairings} applies with \(\dcat = \wcat\) to yield a functor
   \(\vCvcat \times \wCvcat \to \wCvcat\).
  This is given on objects by
  \[
   (\vCvobj, \wCvobj) \mapsto \wCvobj \;\vphantom{\otimes}_R\!\!\otimes_R \vCvobj.
  \]
\end{enumerate}

\subsection{Tall\enhyp{}Wraith Monoids from Free Functors}

The situation studied in \cite{dtgw} was that of commutative, unital rings.
This was extended in \cite{jbbw} to algebras.
Various examples of what we call Tall\enhyp{}Wraith monoids were given in both of those papers.
Of particular relevance are the considerations of group rings.
Indeed, both describe a method of constructing Tall\enhyp{}Wraith monoids from groups.

This is an example of the construction described in Figure~\ref{fig:lifts}.
In Figure~\ref{fig:lifts}, we take a category of commutative, unital algebras for \(\ucat\) and for \(\vcat\) the corresponding category of modules.
Then \(\vCStcat\) is the category of \co{commutative}, \co{unital} coalgebras and \(\vCStStcat\) is the category of unital bialgebras (not necessarily commutative).
Figure~\ref{fig:lifts} tells us that the free commutative algebra on a bialgebra is a Tall\enhyp{}Wraith monoid.

Figure~\ref{fig:lifts} also tells us that the free commutative algebra on a monoid is a Tall\enhyp{}Wraith monoid.
It also tells us that if we start with a group (or monoid), take its free module, resulting in the group ring, and then take the free algebra on that, we obtain the same Tall\enhyp{}Wraith monoid as if we had simply taken the free algebra on our group.

\subsection{Tall\enhyp{}Wraith Monoids from Toy Cohomology Theories}

Section~\ref{sec:toys} gives us another source of relatively simple Tall\enhyp{}Wraith monoids.
As this is closely related to genuine cohomology theories, it is worth taking the time to work through this example carefully.
Let \(\robj\) be a finite \robj.
Let \(\bcat\) be \bcat.
Let us write the Tall\enhyp{}Wraith \(\bcat\)\enhyp{}monoid \(\Hom{\scat}{\robj}{\robj}\) more compactly as \(\robj^{\robj}\).

Now the Tall\enhyp{}Wraith product between \bCbobjalts and \bobjs is defined by the adjunction
\[
  \Hom{\bcat}{\bobj_1 \twprod \bCbobj}{\bobj_2} \cong \Hom{\bcat}{\bobj_1}{\Hom{\bcat}{\bCbobj}{\bobj_2}}
\]
so to understand \(\robj^{\robj} \twprod \bobj\) we can begin by examining \(\Hom{\bcat}{\robj^{\robj}}{\bobj}\), or more generally by examining \(\Hom{\bcat}{\robj^{\sobj}}{\bobj}\) for a finite \sobj, \(\sobj\).

There are several ways to think of \(\Hom{\bcat}{\robj^{\sobj}}{\bobj}\).
A \bmor \(\robj^{\sobj} \to \bobj\) is determined by its effect on the obvious generating set and the image of these generators must satisfy the same relations in the target as they do in the source.
Thus a \bmor, \(\bmor \colon \robj^{\sobj} \to \bobj\), corresponds to an \(\sobj\)\enhyp{}indexed subset of \(\bobj\) of pairwise orthogonal idempotents which sum to \(1\).
Since idempotents in a unital ring correspond to splittings, we can also identify \(\bmor \colon \robj^{\sobj} \to \bobj\) with an isomorphism of rings
\[
  \bobj \cong \prod_{\selt \in \sobj} \bobj_{\selt}.
\]
Of these different views, we shall take the middle one: of (finite) sequences of idempotents.
Thus \(\robj^{\sobj}\) represents the functor that sends a \bobj to the set of \(\sobj\)\enhyp{}indexed sequences of idempotents (that are  orthogonal and sum to \(1\)).

Now if we take \(\sobj = \abs{\robj'}\) where \(\robj'\) is another \robj, we can put the structure of an \(\robj'\)\enhyp{}algebra on the set of \(\abs{\robj'}\)\enhyp{}indexed sequences of idempotents in \(\bobj\).
This is done by gathering suitable terms together:
\begin{align*}
\big((\alpha_{\relt}) + (\beta_{\relt})\big)_{\relt_0} &= \sum_{\relt_1 + \relt_2 = \relt_0} \alpha_{\relt_1} \beta_{\relt_2} \\
\big((\alpha_{\relt}) \cdot (\beta_{\relt})\big)_{\relt_0} &= \sum_{\relt_1 \relt_2 = \relt_0} \alpha_{\relt_1} \beta_{\relt_2} \\
\big(\relt_1 (\alpha_{\relt})\big)_{\relt_0} &= \sum_{\relt_2 : \relt_1 \relt_2 = \relt_0} \alpha_{\relt_2}
\end{align*}
The sequence corresponding to \(\relt \in \robj'\) has the unit in the \(\relt\)\enhyp{}place and zero elsewhere.

As this is a \(\robj'\)\enhyp{}algebra, we can again consider suitable sequences of idempotents in it.
Thus for a \sobj[\sobj] we consider \(\sobj\)\enhyp{}indexed sequences of idempotents of \(\robj'\)\enhyp{}indexed sequences of idempotents in \(\bobj\)!

Given an \(\sobj\)\enhyp{}indexed sequence of idempotents of \(\bobj\), we obtain an \(\sobj\)\enhyp{}indexed sequence of idempotents of \(\robj'\)\enhyp{}indexed sequences of idempotents simply by putting the \(\selt\)\enhyp{}idempotent in \(1\)\enhyp{}component of the sequence corresponding to \(\selt\) and the sum of the rest in the \(0\)\enhyp{}component (so that the sum is \(1\)):
\[
  (\alpha_{\selt}) \mapsto \big((1 - \alpha_{\selt}, \alpha_{\selt}, 0,\dotsc, 0)\big)
\]
where we list \(\robj'\) begining \((0,1,\dotsc)\).

Back in terms of representing objects, this yields a morphism in the other direction:
\[
  \robj^{\robj'} \twprod {\robj'}^{\sobj} \to \robj^{\sobj}.
\]
When \(\robj' = \robj\) and \(\sobj = \abs{\robj}\), we obtain the required morphism
\[
  \mu_{\robj} \colon \robj^{\robj} \twprod \robj^{\robj} \to \robj^{\robj}.
\]

One thing that thinking of \(\robj^{\robj}\) in terms of idempotents brings to the fore is the statement of the following result.

\begin{theorem}
The map
\[
  \mu_{\robj} \colon \robj^{\robj} \twprod \robj^{\robj} \to \robj^{\robj}
\]
is an isomorphism if and only if \(\robj\) has no non\hyp{}trivial idempotents.
\end{theorem}

\begin{proof}
Let \(\bobj\) be a \bobj.
We want to consider \(\robj\)\enhyp{}indexed sequences of idempotents in the ring of \(\robj\)\enhyp{}indexed sequences of idempotents in \(\bobj\).
Thus we need to start with an \(\robj\)\enhyp{}indexed sequence of idempotents in \(\bobj\) and see what it means for that to be idempotent.

So let \((\alpha_{\relt})_{\relt}\) be an \(\robj\)\enhyp{}indexed sequence of idempotents in \(\bobj\).
Recall that we require these to be pairwise orthogonal and to sum to \(1\).
The square of this is given by
\[
  \left( (\alpha_{\relt}) \cdot (\alpha_{\relt}) \right)_{\relt_0} = \sum_{\relt_1 \relt_2 = \relt_0} \alpha_{\relt_1} \alpha_{\relt_2}
\]
For this to be idempotent, we need
\[
  \sum_{\relt_1 \relt_2 = \relt_0} \alpha_{\relt_1} \alpha_{\relt_2} = \alpha_{\relt_0}.
\]

Recall that the sequence \((\alpha_{\relt})_{\relt}\) is itself composed of pairwise orthogonal idempotents.
Thus \(\alpha_{\relt_1} \alpha_{\relt_2}\) is \(\alpha_{\relt_1}\) if \(\relt_1 = \relt_2\) and is \(0\) otherwise.
Hence:
\[
  \left( (\alpha_{\relt}) \cdot (\alpha_{\relt}) \right)_{\relt_0} = \sum_{\relt_1^2 = \relt_0} \alpha_{\relt_1}.
\]
Thus for this sequence to itself be idempotent, we require that for all \(\relt_0 \in \abs{\robj}\),
\begin{equation}
\label{eq:idemofidem}
  \sum_{\relt_1^2 = \relt_0} \alpha_{\relt_1} = \alpha_{\relt_0}.
\end{equation}

Now \(\alpha_{\relt}\) is itself an idempotent and for \(\relt \ne \relt'\), \(\alpha_{\relt} \alpha_{\relt'} = 0\) so multiplying \eqref{eq:idemofidem} by \(\alpha_{\relt_0}\) we see that \(\relt_0\) must appear in the sum on the left\hyp{}hand side.
Thus \(\relt_0^2 = \relt_0\) so \(\relt_0\) is itself an idempotent in \(\robj\).
To show that there are no other terms in the summation in \eqref{eq:idemofidem} we multiply by \(\alpha_{\relt}\) for \(\relt \ne \relt_0\).
On the right, we obtain \(0\); on the left, we obtain \(\alpha_{\relt}\) if \(\relt\) is an index in the summation.
Thus the only index that contributes non\hyp{}trivially is \(\relt_0\) and that only if \(\relt_0\) is an idempotent in \(\robj\).

In particular, if \(\robj\) has no non\hyp{}trivial idempotents then the only idempotent sequences are those with non\hyp{}zero terms in the terms corresponding to \(0\) and to \(1\).
As the sequence must sum to \(1\), the sequence must therefore be of the form
\[
  \big( (1 - \alpha_{\selt}, \alpha_{\selt}, 0, \dotsc, 0)\big)
\]
which is precisely the family of sequences in the image of the map from \(\abs{\robj}\)\enhyp{}indexed sequences.

Hence the map from \(\abs{\robj}\)\enhyp{}indexed sequences of idempotents in \(\bobj\) to  \(\abs{\robj}\)\enhyp{}indexed sequences of idempotents of \(\abs{\robj}\)\enhyp{}indexed sequences of idempotents in \(\bobj\) is a bijection.

For the converse, we assume that \(\robj\) has a non\hyp{}trivial idempotent, say \(\alpha\).
Order \(\abs{\robj}\) starting with \((0,1,\alpha, 1 - \alpha, \dotsc)\).
Then consider the sequence of idempotents in \(\robj\):
\[
  (0,0,\alpha, 1 - \alpha, 0, \dotsc).
\]
These are pairwise orthogonal and sum to \(1\).
Squaring this produces an idempotent.
Together with
\[
  (0,0,1 - \alpha, \alpha, 0, \dotsc)
\]
we get a pair of idempotent sequences which are orthogonal and sum to \(1\).
These two are not in the image of the map from \(\abs{\robj}\)\enhyp{}indexed sequences, and hence that map is not a bijection.
\end{proof}

Also included in the Tall\enhyp{}Wraith structure is a map from the object representing the identity functor.
For \bobjs, this is \(\robj \lb\selt\rb\), the free \(\robj\)\enhyp{}algebra on a single generator.
The morphism \(\eta_{\robj} \colon \robj\lb\selt\rb \to \robj^{\robj}\) is simple to describe: it is determined by sending  \(\selt\) to the identity morphism \(\robj \to \robj\).

This morphism is exactly the type of situation that we wish to study.
We have an unknown thing, \(\robj^{\robj}\), and a known thing, \(\robj\lb\selt\rb\), and a morphism from the known to the unknown.
The questions then become: ``How much of \(\robj^{\robj}\) does this morphism see?'' and ``What of \(\robj\lb\selt\rb\) is redundant information?''.
In other words, what are the image and kernel of \(\eta_{\robj} \colon \robj\lb\selt\rb \to \robj^{\robj}\).

The situation is particularly simple in the case of finite fields and we end this paper with a nice description of the Tall\enhyp{}Wraith monoid \(\Hom{\scat}{\robj}{\robj}\) in this case.

So let \(\robj\) be a finite field.
Let \(\robj\) have characteristic \(p\) and order \(p^n\).
As the non\hyp{}zero elements of \(\robj\) form a cyclic group of order \(\abs{\robj} - 1\), we see that the morphism \(\robj \to \robj\) corresponding to the image of \(\selt^{p^n - 1}\) takes \(0\) to itself and everything else to \(1\).
Therefore \(\delta_0 = 1 - \selt^{p^n - 1}\) and once we have \(\delta_0\), we get everything else by pre\hyp{} or post\hyp{}composition.

The kernel is the ideal generated by \(\selt^{p^n} - \selt\).
Let us show directly that the quotient \(\robj\lb\selt\rb/(\selt^{p^n} - \selt)\) inherits a Tall\enhyp{}Wraith monoid structure from \(\robj\lb\selt\rb\).

Firstly, by construction \((\selt^{p^n} - \selt)\) is an ideal, so the quotient inherits a ring structure.

Secondly, we show that the \co{ring} structure passes to the quotient.
To do this, we first compute the \co{multiplication} and \co{addition}.
This uses the fact that both are morphisms of \(\robj\)\enhyp{}algebras and that \(\selt\) is \emph{ring\enhyp{}like}.
Thus we obtain the following formul\ae.

\begin{align*}
\Delta^+ (\selt^{p^n} - \selt) &= \left(\Delta^+ \selt\right)^{p^n} - \Delta^+ \selt \\
&= \left( 1 \otimes \selt + \selt \otimes 1 \right)^{p^n} - 1 \otimes \selt - \selt \otimes 1 \\
&= 1 \otimes \selt^{p^n} + \selt^{p^n} \otimes 1  - 1 \otimes \selt - \selt \otimes 1 \\
&= 1 \otimes (\selt^{p^n} - \selt) + (\selt^{p^n} - \selt) \otimes 1 \\
\Delta^\times (\selt^{p^n} - \selt) &= \left(\Delta^\times \selt\right)^{p^n} - \Delta^\times \selt \\
&= \left( \selt \otimes \selt \right)^{p^n} - \selt \otimes \selt \\
&= \selt^{p^n} \otimes \selt^{p^n} - \selt \otimes \selt \\
&= \selt^{p^n} \otimes \selt^{p^n} - \selt^{p^n} \otimes \selt + \selt \otimes \selt^{p^n} - \selt \otimes \selt \\
&= \selt^{p^n} \otimes( \selt^{p^n} - \selt) + \selt \otimes (\selt^{p^n} - \selt)
\end{align*}

For \(\relt \in \robj\), the induced morphism \(\epsilon_{\relt} \colon \robj\lb\selt\rb \to \robj\) is the evaluation map which clearly takes \(\selt^{p^n} - \selt\) to \(0\), regardless of \(\relt\).
Note that this also takes care of the \co{units} for \co{addition} and \co{multiplication} as these are \(\epsilon_0\) and \(\epsilon_1\) respectively.
It also takes care of the \co{inverse} for \co{addition} as that is \co{multiplication} by \(\epsilon_{-1}\).

Therefore the ideal generated by \(\selt^{p^n} - \selt\) is a \co{ideal} for both \(\Delta^+\) and \(\Delta^\times\), hence the quotient inherits the \co{operations}.

Finally, we consider the Tall\enhyp{}Wraith product.
Given a Tall\enhyp{}Wraith monoid in \bcat, say \(\bCbTtobj\), and an ideal \(I\) which is also a \co{ideal}, for the Tall\enhyp{}Wraith product to descend to the quotient, we must have that the following holds for all \(\bCbTtelt_1 \in \bCbTtobj\) and \(\bCbTtelt_2 \in I\):
\[
  \bCbTtelt_1 \twprod \bCbTtelt_2 - \bCbTtelt_1 \twprod 0 \in I.
\]

In \(\robj[x]\), the Tall\enhyp{}Wraith product is composition of polynomials.
Since, in this case, \(I\) is the ideal generated by \(\selt^{p^n} - \selt\), it contains all powers and multiples of \(\selt^{p^n} - \selt\).
Thus if \(q(\selt)\) is a polynomial with no constant term, \(q(\selt^{p^n} - \selt) \in I\).
More generally, for any polynomial, \(q(\selt^{p^n} - \selt) - q(0) \in I\), as required.

Thus we have proved the following special case of the first isomorphism theorem for Tall\enhyp{}Wraith monoids.

\begin{theorem}
If \(\robj\) is a finite field, then there is an isomorphism of Tall\enhyp{}Wraith monoids
\[
  \eta_{\robj} \colon \robj \lb\selt\rb/(\selt^{p^n} - \selt) \to \Hom{\scat}{\abs{\robj}}{\abs{\robj}},
\]
which sends (the equivalence class of) \(\selt\) to the identity morphism on \(\abs{\robj}\). \noproof
\end{theorem}

\appendix

\section{Notation}
\label{ap:notation}

This paper is replete with abstract categories that exist in various relationships to each other.
We have tried to find a notation that is both concise and conveys the relationships between the various categories.
As no such system is going to be perfect, and as what may seem clear to us may not be clear to the casual reader, we include a brief explanation of our conventions here as a reference.
This notation is the same as that used in our previous article, \cite{MR2559638}, although in that article there were even more variations.

We work with categories that are built out of other categories by taking certain objects with extra structure.
Usually, this extra structure is encoded in a (finite) algebraic theory, which can itself be viewed as a category.
Thus at the simplest level, we begin with a category \(\dcat\) and an algebraic theory \(\vcat\) and form the category \dvcat which we write as \(\dvcat\).

We frequently need to iterate this construction, but with modifications.
For example, to consider \Cvalgobjs we need to first take the opposite category of the base category (which we denote by \(\Odcat\)), take \valgobjs, and then take the opposite category again.
We write the resulting category as \(\dCvcat\).
When iterating these constructions, we read from right to left.
Thus \(\dvCvcat\) is \dvCvcat.

Another modification we use is when we have a monoidal category.
Then, assuming that our algebraic theory comes from an operad, we can work with the monoidal product rather than the cartesian product.
To denote this modification, we use a decoration: \(\mTo{}\) (or \(\mSo{}\) depending on the notation for the monoidal product) and so write \(\dTtcat\) for the category of \talgobjs in \dcat using the \(\twprod\)\enhyp{}product instead of the cartesian product.

The last decorations that we use in this paper relate to graded and ungraded objects.
As it is important to keep in mind when something is graded or ungraded, we use unadorned notation \(\vcat\) to refer to ungraded algebraic theories and denote graded ones by \(\mGo{}\), thus \(\Gvcat\) is a graded theory.
Note that we do not assume a connection between \(\vcat\) and \(\Gvcat\) (indeed, we never use the same letter for both in the same context).
Where we do want to take a category and impose a grading upon it, we use the notation \(\mZo{}\).
Thus \(\Zdcat\) denotes \Zdcat.

Most of our categories are very general.
A few have special meaning.
These, together with a list of the above modifications, are gathered into Table~\ref{tab:notation}.

\begin{table}
\begin{centre}
\begin{tabular}{rl}
\toprule
Symbol & Meaning \\
\midrule
\(\dcat\) & Base category \\
\(\Odcat\) & \Odcatu \\
\(\Zdcat\) & \Zdcatu \\ 
\(\dvcat\) & \dvcatu \\
\(\dCvcat\) & \dCvcatu \\
\(\dGvcat\) & \dGvcatu \\
\(\dCGvcat\) & \dCGvcatu \\
\(\vCvTtcat\) & \vCvTtcatu \\
\(\rcat\) & \rcatu \\
\(\bcat\) & \bcatu \\
\(\tcat\) & The category of monoids (in \scat)\\
\(\scat\) & \scatu \\
\bottomrule
\end{tabular}
\end{centre}
\caption{Table of notatation}
\label{tab:notation}
\end{table}

\bibliography{bibtex_export}

\begin{thebibliography}{KPMS82}

\bibitem[Ber98]{gb}
George~M. Bergman.
\newblock An invitation to general algebra and universal constructions.
\newblock Berkeley, CA, 1998.

\bibitem[BH96]{gbah}
George~M. Bergman and Adam~O. Hausknecht.
\newblock Co-groups and co-rings in categories of associative rings.
\newblock Providence, RI, 1996.

\bibitem[BW05]{jbbw}
James Borger and Ben Wieland.
\newblock Plethystic algebra.
\newblock {\em Adv. Math.}, 194:246–283, 2005.

\bibitem[Fre66]{pf2}
P.~Freyd.
\newblock Algebra valued functors in general and tensor products in particular.
\newblock {\em Colloq. Math.}, 14:89–106, 1966.

\bibitem[Kan58]{MR0111035}
Daniel~M. Kan.
\newblock On monoids and their dual.
\newblock {\em Bol. Soc. Mat. Mexicana (2)}, 3:52–61, 1958.

\bibitem[KPMS82]{wkjmmpis}
W.~Kühnel, M.~Pfender, J.~Meseguer, and I.~Sols.
\newblock Algebras with actions and automata.
\newblock {\em Internat. J. Math. Math. Sci.}, 5:61–85, 1982.

\bibitem[Lin66]{fl2}
F.~E.~J. Linton.
\newblock Autonomous equational categories.
\newblock {\em J. Math. Mech.}, 15:637–642, 1966.

\bibitem[SW08]{assw4}
Andrew Stacey and Sarah Whitehouse.
\newblock Stable and unstable operations in mod {$p$} cohomology theories.
\newblock {\em Algebr. Geom. Topol.}, 8:1059–1091, 2008.

\bibitem[SW09]{MR2559638}
Andrew Stacey and Sarah Whitehouse.
\newblock The hunting of the {H}opf ring.
\newblock {\em Homology, Homotopy Appl.}, 11:75–132, 2009.

\bibitem[TW70]{dtgw}
D.~O. Tall and G.~C. Wraith.
\newblock Representable functors and operations on rings.
\newblock {\em Proc. London Math. Soc. (3)}, 20:619–643, 1970.

\end{thebibliography}

\end{document}